\documentclass[11pt,leqno]{article}
\def\s{\dot{s}}
\def\w{\dot{w}}
\def\v{{\rm v}}
\def\sigmad{\dot{\sigma}}
\def\int{\mathbb{Z}}

\def\Ue{{\cal U}_{\varepsilon}({\mathfrak g})}

\def\O{{\cal O}}

\def\proof{{\bf Proof. }}

\def\Pf{\proof}
\def\Proof{\proof}
\def\pf{\proof}
\def\id{{\rm id}}
\def\V{{\cal V}}

\usepackage{times}
\usepackage{graphics}
\usepackage{amssymb}
\usepackage{amscd}
\usepackage{amsmath}
\input xy
\xyoption{all}
\title{Spherical conjugacy classes and involutions in the Weyl group}

\newtheorem{theorem}{Theorem}[section]
\newtheorem{lemma}[theorem]{Lemma}
\newtheorem{corollary}[theorem]{Corollary}

\newtheorem{definition}[theorem]{Definition}
\newtheorem{remark}[theorem]{Remark}
\newtheorem{example}[theorem]{Example}

\author{Giovanna Carnovale\\
Dipartimento di Matematica Pura ed Applicata\\
Torre Archimede - via Trieste 63 - 35121 Padova - Italy\\
email: carnoval@math.unipd.it }
\date{}
\begin{document}
\maketitle
\begin{abstract}Let $G$ be a simple algebraic group over an
  algebraically closed field of characteristic zero or positive odd, good characteristic. Let
  $B$ be a Borel subgroup of $G$. We show that the spherical conjugacy
  classes of $G$  intersect only the double cosets of $B$ in $G$
  corresponding to involutions in the Weyl group of $G$. This result
  is used in order to prove that for a spherical conjugacy class $\O$
  with dense $B$-orbit $v_0\subset BwB$ there holds $\ell(w)+{\rm
  rk}(1-w)=\dim\O$ extending a characterization of spherical conjugacy
  classes obtained by Cantarini, Costantini and the author to
  the case of groups over fields of odd, good characteristic.  
\end{abstract}

\section*{Introduction}

If an algebraic group acts with finitely many orbits, a natural way to
 understand the action is given by the combinatorics of the Zariski
 closures of such orbits. In \cite{RS}, \cite{results}, a detailed description of the
 combinatorics of the closures of orbits for a Borel
 subgroup $B$ in a symmetric space $G/K$ is
 given. The description is provided in terms of an action, on the
 set of these orbits, of a monoid $M(W)$ related to the Weyl group $W$
 of $G$. This action  
 is best understood considering the decomposition into $B$-orbits of an orbit of
 a minimal parabolic subgroup. Through this approach 
 several invariants of the $B$-orbits can be determined, including their
 dimension. To each $B$-orbit it is possible to associate a Weyl group
 element and  the Weyl group
 element corresponding to the (unique) dense $B$-orbit in the
 symmetric space can be described in combinatorial terms. 
 A formula for the dimension of each $B$-orbit is provided in terms of its associated Weyl
 group element and the sequence of elements in the monoid that are
 necessary to reach it from a closed $B$-orbit. When the symmetric space corresponds to an
 inner involution, that is, if it corresponds to a conjugacy
 class in $G$, the attached Weyl group element is just the 
element corresponding to the Bruhat cell containing the $B$-orbit. 

The monoid action can be carried over to homogeneous spaces of
algebraic groups for which
the action of the Borel subgroup has finitely many orbits, i.e., the
spherical homogeneous spaces (\cite{montreal}) and it can be
used to define representations of the Hecke algebra (\cite{MS}). A
more geometric approach to a Bruhat order on spherical varieties has
been addressed in \cite{B-bruhat}. Besides, a genuine Weyl group
action on the set of $B$-orbits on a spherical homogeneous space was
defined in \cite{knop}.

The action of $M(W)$ on a
spherical homogeneous space does not afford all nice properties that it
had in the symmetric case (see \cite{closures} for a few
key counterexamples) and it is natural to ask which properties still
hold for spherical conjugacy classes. One of the main differences
between the geneal spherical case and the symmetric case is that there
are $B$-orbits that do not lie in the orbit of a closed one when
acting by $M(W)$. However, every $B$-orbit can be reached from a
closed $B$-orbit through a sequence of moves involving either the
$M(W)$-action or the $W$-action (\cite{montreal}).
 
A natural question is whether we can provide formulas for the
dimension of each $B$-orbit in a spherical conjugacy class 
in terms of the actions of $M(W)$ and
$W$. Although not all results in \cite{RS} 
hold at this level of generality, there are properties that
hold true in general. For instance, the dimension of the dense
$B$-orbit in a spherical conjugacy class is governed by a formula
analogous to the formula for the dimension of a $B$-orbit in a
symmetric conjugacy class. This result was achieved in \cite{ccc}, leading to a
characterization of spherical conjugacy classes in simple algebraic
groups. The interest in this formula lied in the 
verification of De Concini-Kac-Procesi conjecture on the
dimension of irreducible representations of quantum groups at the
roots of unity (\cite{DC-K-P}) in the case of spherical conjugacy
classes. For this reason, the analysis was restriced to the case of an
algebraic group over an algebraically closed field of characteristic
zero. In order to obtain the characterization, a classification of all
spherical conjugacy classes in a simple algebraic group was needed,
and part of the results were obtained through a case-by-case analysis
involving this classification. 
 
In the present paper we apply the combinatorics of $B$-orbit closures
for spherical conjugacy classes in simple algebraic groups  to retrieve 
the formula in \cite{ccc}. This will show that the characterization of
spherical conjugacy classes can be
achieved without
making use of their classification and
without drastic restrictions on the characteristic of the base field. 

A first question to be answered concerns which Weyl group elements may
correspond to a $B$-orbit in a spherical conjugacy class. In the case
of a symmetric conjugacy class it is immediate to see that such Weyl group
elements are involutions. An analysis of the actions of
$M(W)$ and $W$ allows us to generalize this result to all spherical
conjugacy classes.

\smallskip

\noindent{\bf Theorem 1 } {\em All $B$-orbits in a spherical conjugacy
  class lie in Bruhat cells corresponding to involutions in the Weyl group}. 

\smallskip

In order to understand the Weyl group elements associated with the
dense $B$-orbit wee analyze the variation of the Weyl group element 
with respect to the action of the monoid $M(W)$. This
analysis leads to a description of the stationary points, i.e., of
those $B$-orbits for which the associated Weyl group element does not
change under the action of all standard generators of $M(W)$. 
Stationary point different from the dense $B$-orbit cannot exist
in symmetric conjugacy classes but they exist, for instance, in
spherical unipotent conjugacy classes. 

The results in
\cite{results} allow us to describe the Weyl group
element corresponding to a stationary point, and more precisely, 
the one associated with the dense $B$-orbit.

Combining the analysis of representatives of the dense $B$-orbit with
results in \cite{closures} yields a new proof of the
formula in \cite{ccc}, that holds now in almost all characteristcs and
does not require the classification of spherical conjugacy classes:

\smallskip

\noindent{\bf Theorem 2} {\em Let $\O$ be a spherical conjugacy class
  in a simple algebraic group $G$,  let $v_0$ be its dense $B$-orbit
  and let $BwB$ be the Bruhat double coset containing $v_0$. Then
  $\dim\O=\ell(w)+{\rm rk}(1-w)$. } 

\smallskip

It is proved in \cite{ccc} with a characteristic-free argument that if
a conjugacy class $\O$ intersects some $BwB$  with $\ell(w)+{\rm
  rk}(1-w)=\dim\O$ then $\O$ is spherical, hence the results in the
present paper provide a characteristic-free proof of the characterization of
spherical conjugacy class given in \cite{ccc}. 

\section*{Acknowledgements}

This research was partially supported by Marie Curie Research Training Network
LIEGRITS Flags, Quivers and Invariant Theory in Lie Representation Theory
number MRTN-CT 2003-505078. 

\section{Preliminaries}\label{preliminaries}

Let $G$ be a simple algebraic group over an algebraically closed field
$k$ of characteristic $0$ or odd and good (\cite[\S 4.3]{131}). Let $B$ be a
Borel subgroup of $G$, $T$ a maximal torus contained in $B$ and $B^-$ the
Borel subgroup opposite to $B$. Let $U$
(respectively $U^-$) be the unipotent radical of $B$ (respectively
$B^-$). For a group $K$ we shall denote by $K^\circ$ the connected
component of the identity of $K$.

We shall denote by  $\Phi$ the  set of roots
relative to $(B, T)$, by $\Phi^+$ the corresponding positive roots, by
$\Delta=\{\alpha_1, \dots, \alpha_n\}$
the corresponding set of simple roots. We shall use the numbering of
the simple roots in \cite[Planches I-IX]{bourbaki}. The usual inner product
between the roots $\alpha$ and $\beta$ will be denoted by
$(\alpha,\beta)$. If the root system is simply laced we will consider
all its roots as long roots. We shall indicate by $P^+$ and $Q$,
respectively, the set of dominant weights and the root lattice
associated with $\Phi$. An element of $T$ corresponding to the co-root $\alpha^\vee$
will be indicated by $\alpha^\vee(h)$.

By $\Pi$ we shall always denote
a subset of $\Delta$ and $\Phi(\Pi)$ will indicate the corresponding
root subsystem.

We shall denote by $W$ be the
Weyl group associated with $G$ and by
$s_{\alpha}$ the reflection corresponding to the root $\alpha$. By
$\ell(w)$ we shall denote the length of the element $w\in W$ and
by ${\rm rk}(1-w)$ we shall mean the rank of $1-w$ in the
standard representation of the Weyl group. By $w_0$ we shall denote
the  longest element in $W$ and by $\vartheta$ we shall denote the
automorphism of $\Phi$ given by $-w_0$. If $\Pi\subset\Delta$ we shall denote by
$W_\Pi$ the parabolic subgroup of $W$ generated by the simple
reflections in $\Pi$. If $N(T)$ is the normalizer
of $T$ in $G$  then $W=N(T)/T$; given an element $w\in W$ we shall denote
a representative of $w$ in $N(T)$ by $\w$.
For any root $\alpha$ of $\Phi$ we shall denote by $x_{\alpha}(t)$ the
elements of the corresponding root subgroup $X_{\alpha}$ of $G$. We
shall choose the representatives $\dot{s}_{\alpha}\in N$ of the reflection
$s_{\alpha}\in W$ as in \cite[Theorem 8.1.4]{springer}.

We assume that we have fixed an ordering of the positive roots so that
every $u\in U$ is written uniquely as an ordered product of element of
the form $x_\alpha(t)$, for $t\in k$ and $\alpha\in\Phi^+$. Given an
element $u\in U$ (respectively $U^-$), by abuse of language we will say that
a root $\gamma\in \Phi^+$ (respectively $-\Phi^+$) {\em occurs in $u$} if
for the expression of $u$ as an ordered product of
$x_\alpha(t_\alpha)$'s  we have $t_\gamma\neq0$. 

If $\alpha\in\Delta$ we shall indicate by $P_\alpha$ the minimal non
solvable parabolic subgroup containing $X_{-\alpha}$ and by
$P_\alpha^u$ its unipotent radical.

For $w\in W$, we shall denote by $U^w$ (respectively, $U_w$) the
subgroup generated by the root subgroups $X_\alpha$  corresponding to
those $\alpha\in\Phi^+$ for which $w^{-1}(\alpha)\in -\Phi^+$
(respectively, $\Phi^+$). We shall denote by $T^w$ the subgroup of the
torus that is centralized by any representative $\w$ of $w$. 

Given an element $x\in G$ we shall denote by $\O_x$ the conjugacy
class of $x$ in $G$ and by $G_x$
(resp. $B_x$, resp. $T_x$) the centralizer of $x$ in $G$ (resp. $B$, resp. $T$). For a
homogeneous space $G/H$ we shall denote by $\V$ the set of $B$-orbits
in $G/H$.  

\begin{definition}\label{sferica} Let $K$ be a connected algebraic group over $k$ and let $H$ be a closed subgroup of $K$.
The homogeneous space $K/H$ is called spherical if there exists a Borel subgroup of $K$ with  a dense orbit.
\end{definition}

It is well-known (\cite{Bri}, \cite{Vin} in characteristic $0$,
\cite{gross}, \cite{knop} in positive characteristic)
that, if $G/H$ is a spherical homogeneous space the set $\V$ of
$B$-orbits in $G/H$ is finite. 

\section{$B$-orbits and Bruhat decomposition}

Let $\O$ be a conjugacy class of $G$ and let $\V$ be the set
of $B$-orbits in $\O$. There is a natural map 
$\phi\colon \V\to W$ associating to 
$v\in\V$ the element $w$ in the Weyl group of $G$ for which
$v\subset BwB$.  The set $\V$ carries a partial order given by:
$v\leq v'\quad{\mbox{ if }}\quad
\overline{v}\subset\overline{v'}$.
When $\O$ is spherical the minimal $B$-orbits are the closed ones and
there is a unique maximal orbit, namely the dense $B$-orbit $v_0$ in $\O$. 

\begin{lemma}\label{le}Let $\O$ be a conjugacy class and let $v,v'\in
  \V$. If $v\leq v'$ then $\phi(v)\leq \phi(v')$ in the
 Bruhat order in $W$. 
\end{lemma}
\Pf We have:
$v\subset \overline{v}\subset\overline{v'}\subset\overline{B\phi(v')B}=
\dot{\cup}_{\sigma\le\phi(v')}B\sigma B$
so $\phi(v)\leq\phi(v')$.\hfill$\Box$

\smallskip 

\begin{lemma}\label{uno} Let $x\in G$ be either semisimple or unipotent and let $\O_x$ be a spherical 
  conjugacy class. The image through $\phi$ of a closed $B$-orbit in $\O$ is $1$. 
\end{lemma}
\Pf By \cite[\S 3.4 (b)]{montreal} if $\O=\O_x$ and $H=G_x$, 
the closed $B$-orbits correspond to those $BgH/H$ with $BgH$ closed in
$G$, so that $B_H=(H\cap g^{-1}Bg)^\circ$ is a Borel subgroup of $H$ and
of $H^\circ$. Let $x$ be semisimple. Since it is not restrictive to
assume that $G$ is simply connected, we have $H=H^\circ$ and $x\in
Z(H)=Z(B_H)$ by \cite[Corollary 6.2.9]{springer}. Hence, the representative $gxg^{-1}$ of the closed
$B$-orbit lies in $B$. 

Let $x$ be unipotent. By \cite[\S 3.15]{131} with $S=\{x\}$ the element $x\in H^\circ$, hence $x\in
Z(H^\circ)\subset B_H$ and we have the statement.\hfill$\Box$

\smallskip

\begin{remark}{\rm All closed $B$-orbits in a spherical conjugacy class $\O_x$ have the
    same dimension (\cite[\S 3.4 (b)]{montreal}) namely $\dim B-\dim B_{G_x}$ where $B_{G_x}$ denotes a
    Borel subgroup in the centralizer of $x$.} 
\end{remark}

\begin{remark}{\rm The converse of Lemma \ref{uno} does not hold for
    spherical unipotent elements. For
    instance, if $\O$ is a spherical unipotent conjugacy class in
    $G=SL_n({\mathbb C})$ the combinatorics of the closures of the
    $B$-orbits that are contained in $B$ is described in
    \cite{anna1}: if $\O$ corresponds to the minimal nilpotent orbit
    in $G$ the $B$-orbits that are contained in $B$ are in bijection
    with the traspositions in $S_n$, and only one, namely
    $B.x_{\beta_1}(1)=X_{\beta_1}\setminus\{1\}$, is closed.} 
\end{remark}

\begin{remark}{\rm In a spherical semisimple conjugacy class, 
$v$ being closed is equivalent to $\phi(B.x)=1$ and to $v\cap T\neq\emptyset$. 
Indeed if $v$ is closed then $v\subset B$ so a representative $x\in v$
is conjugate in $B$ 
to some element in $T$. Viceversa, if $v=B.t$ for some $t\in
T$ then $t$ normalizes $B$ and by \cite[Theorem 9.2]{borel} the
conjugacy class $B.t$ is closed.}
\end{remark}

\smallskip 

Let $M=M(W)$ be the monoid with elements $m(w)$ indexed by the
elements $w\in W$ with relations
$$m(s)m(w)=m(sw), \mbox{ if } l(sw)>l(w),\;m(s)m(w)=m(w), \mbox{ if }
l(sw)<l(w). $$ 
The monoid $M(W)$ is generated by the elements $m(s)$ corresponding to
simple reflections, subject to the braid relations and
to the relation $m(s)^2=m(s)$. In \cite{RS} an action of the monoid
$M(W)$ on the set of
$B$-orbits of a symmetric space $G/K$ is defined. This action can be generalized
to an action of $M(W)$ on the set $\V$ of $B$-orbits of a spherical
homogeneous space (see, for instance, \cite[\S
  3.6]{montreal}). The action of $m(s)$, for a
simple reflection $s=s_\alpha$ is given as follows. If $P_\alpha$
is the minimal parabolic subgroup corresponding to $\alpha$ and $v\in
\V$ then
$m(s).v$ is the dense $B$-orbit in $P_\alpha v$. This action
is analyzed in \cite{closures}, \cite{knop},\cite[\S 4.1]{MS}, \cite{RS}.
We provide an account of the information we will need.

Given $v\in \V$, choose $y\in v$ with stabilizer $(P_\alpha)_y$ in
$P_\alpha$. Then $(P_\alpha)_y$ acts on $P_\alpha/B\cong {\mathbb P}^1$ with finitely
many orbits. Let $\psi\colon (P_\alpha)_y\to PGL_2(k)$ be the corresponding
group morphism. The kernel of $\psi$ is ${\rm Ker}(\alpha)P^u_\alpha$.


The image $H$ of $(P_\alpha)_y$ in $PGL_2(k)$ is either:
$PGL_2(k)$; or solvable and contains a nontrivial unipotent subgroup;
or a torus; or the normalizer of a torus. Here is a list of the possibilities
that may occur:

\begin{enumerate}
\item[I] $P_\alpha v=v$ so $H=PGL_2(k)$; 
\item[IIa]$P_\alpha v=v\cup m(s)v$, with $\dim v=\dim P_\alpha
  v-1$. In this case we may choose an $y\in v$ for which $\psi(X_\alpha)\subset
  H\subset\psi(B)$. 
\item[IIb] $P_\alpha v=v\cup v'$, with $\dim v'=\dim v-1$ and $v$ is open in
  $P_\alpha v$ so $m(s)v=v$. In this case
  we may choose an $y\in v$ for which $\psi(X_{-\alpha})\subset
  H\subset\psi(B^-)$.  
\item[IIIa] $P_\alpha v=v\cup v'\cup m(s)v$, with $\dim v=\dim v'=\dim
  P_\alpha v-1$
  and $v\neq v'$. In this case
  we may choose an $y\in v$ for which $H=\psi(T)$. Then
  $\s_\alpha y\s_{\alpha}^{-1}\in v'$ and $\s_\alpha x_\alpha(t) y
  x_{\alpha}(-t)\s_{\alpha}^{-1}\in m(s)v$ for $t\neq0$.
\item[IIIb] $Pv=v\cup v'\cup v''$, with $\dim v-1=\dim v'=\dim v''$
  and $v$ is open in
  $P_\alpha v$ so $m(s)v=v$. In this case
  we may choose an $y\in v$ for which $H=\psi(\s_\alpha x_{\alpha}(-1)
  Tx_{\alpha}(1)\s_\alpha^{-1})$. 
\item[IVa] $Pv=v\cup m(s)v$, with $\dim v=\dim P_\alpha v-1$. We may
 choose an $y\in v$ for which $H=\psi(N(T))$.  
 Besides, $\s_\alpha y\s_{\alpha}^{-1}\in v$ and $\s_\alpha x_\alpha(t) y
  x_{\alpha}(-t)\s_{\alpha}^{-1}\in m(s)v$ for $t\neq0$.
\item[IVb] $Pv=v\cup v'$, with $\dim v=\dim v' +1$, and $v$ is open in
  $P_\alpha v$ so $m(s)v=v$.  In this case
  we may choose an $y\in v$ for which $H=\psi(N(\s_\alpha
  T\s_\alpha^{-1}))$. 
\end{enumerate}
Based on the structure of $H$, cases II, III, and IV are also called
type $U$, type $T$ and type $N$, respectively.

Based on this analysis a $W$-action on $\V$ can be defined (\cite{knop}, \cite[\S 4.2.5,
  Remark]{MS}) as follows: in case II the two $B$-orbits are
interchanged; in case III the two non-dense orbits are interchanged,
in all other cases the $B$-orbits are fixed. The image of $v\in\V$
  through the action of a simple reflection $s\in W$ will be denoted
  by $s.v$. 

We recall (\cite[\S 3.6]{montreal}) that a {\em reduced decomposition}
of $v\in \V$ is a pair $({\bf v},{\bf s})$ with ${\bf v}=(v(0),v(1),\ldots, v(r))$ a sequence of
distinct elements in $\V$ and ${\bf s}=(s_{i_1},\ldots, s_{i_r})$ a
sequence of simple reflections such that $v(0)$ is closed and
$v(j)=m(s_{i_j}).(v({j-1}))$ for $1\leq j\leq r-1$.

All $B$-orbits in a symmetric homogeneous space admit a reduced
decomposition (\cite[\S 7]{RS}). This is still the case for the dense
$B$-orbit in spherical homogeneous spaces but it is not always the case 
for generic $B$-orbits in a spherical homogeneous space. The reader can refer
to \cite{closures} for a series of counterexamples. We will use a
weaker notion of decomposition that exists for every $v\in \V$. 

Given a reduced decomposition $({\bf v},{\bf s})=((v(0),\ldots, v(r)), {\bf s})$ of $v\in \V$ a {\em subexpression} of $({\bf v},{\bf s})$ 
(\cite[\S 3.6]{montreal}) 
 is a sequence 
${\bf  x}=(v'(0),v'(1),\ldots,\,v'(r))$ of elements in $\V$ with
 $v'(0)=v(0)$ and such that for $1\leq
i\leq r$ only one of the following alternatives occurs:
\begin{enumerate}
\item[(a)] $v'(j-1)=v'(j)$; 
\item[(b)] $v'(j-1)\not=v'(j)$, $\dim v'(j-1)=v'(j)$ and $v'(j)=s_j.v'(j-1)$;
\item[(c)] $\dim v'(j-1)=\dim v'(j)-1$ and $v'(j)=m(s_j).(v'(j-1))$. 
\end{enumerate}
The element $v'(r)$ is called the {\em final term} of the
subexpression. Even though some $B$-orbits in a spherical homogeneous
space might not have a reduced decomposition, every $v\in \V$ is the
final term of a subexpression of a reduced decomposition of the
dense $B$-orbit $v_0$. This is to be found in \cite[\S 3.6 Proposition
  2]{montreal} and it holds also in positive odd characteristic.

\begin{lemma}\label{invo}Let $\O$ be a spherical conjugacy class, let $v\in \V$ and let $s$
  be a simple reflection. If $w=\phi(v)$ is an involution then
  $w'=\phi(m(s).v)$ is an involution. 
\end{lemma}
\Pf We consider $P_\alpha v=v\cup
(BsB).v$. If $m(s)v=v$ there is nothing to prove. 
Let us assume that $m(s)v\subset
(Bs).v$. Then $m(s).v\subset BsBwBsB$. The following four
possibilities may occur:  
\begin{enumerate}
\item $l(sws)=l(w)+2$: by \cite[Lemma 8.3.7]{springer} we have
$w'=sws$ and the statement holds. 
\item $l(sw)>l(w)$ and $l(sws)=l(w)$: since $w$ is an involution, by \cite[Lemma 3.2 (ii)]{results} for $\theta=\id$ 
we have $sw=ws$. Besides, by \cite[Lemma 8.3.7]{springer} we have 
$$m(s).v\subset BswBsB=BswB\cup BswsB.$$ hence $w'\in\{sw,w\}$ is an involution.
\item $l(sw)<l(w)$ and $l(sws)=l(w)$: again by  \cite[Lemma 3.2 (i)]{results} for $\theta=\id$ 
we have $sw=ws$. Then $$m(s).v\subset BsBwBsB\subset
BswsB\cup BswB$$ so $w'\in \{w, ws\}$ is an involution.
\item $l(sws)=l(w)-2$: we have 
$$m(s).v\subset BsBwBsB\subset BswB\cup BswsB\cup Bws B\cup BwB.$$  
By \cite[\S 3.6 Proposition 1 (a)]{montreal} we have $v\leq
m(s).v$ so by Lemma \ref{le} there holds $w\leq w'$ hence $w'\neq sw, ws,
sws$ and $w'=w$ is an involution.\hfill$\Box$ 
\end{enumerate}

\begin{corollary}\label{reduced}Let $\O_x$ be a spherical
  conjugacy class  with $x$ either semisimple or unipotent and let
  $v\in \V$ admit a reduced
  decomposition. Then $\phi(v)$ is an involution. In particular, this
  holds for the dense $B$-orbit of a spherical semisimple or of a
  spherical unipotent conjugacy class. 
\end{corollary}
\Pf The first assertion follows by induction on $\dim v$ using Lemma
\ref{uno} and Lemma \ref{invo}. 
The last statement follows from \cite[\S 3.6, Proposition
  3(ii)]{montreal}.\hfill$\Box$ 

\smallskip

\begin{theorem}\label{image}Let $\O$ be a spherical conjugacy class,
 and let $\phi\colon \V\to W$ be the natural
  map. Then the image of $\phi$ consists of involutions. 
\end{theorem}  
\Pf We first consider a spherical semisimple conjugacy class. Let 
$v\in \V$ and let ${\bf x}$ be a subexpression of a reduced
decomposition of the dense $B$-orbit $v_0$ with initial term a closed $B$-orbit
$v(0)$  and final term $v$. 
We proceed by induction on
$\dim v$. If $v$ has minimal dimension then it is closed, otherwise
it would contain in its closure a $B$-orbit of strictly smaller
dimension. It follows from Lemma \ref{uno} that $\phi(v)=1$. 

Let $v=v(r)$ not be closed so that $\dim v(r)>\dim v(0)$.  If
$v(r)=v(r-1)$ we may shorten the subexpression replacing $r$ by
$r-1$. Hence we may assume that $v(r)=m(s)v(r-1)$ or $v(r)=s.v(r-1)$
for some simple reflection $s$. If $v(r)=m(s)v(r-1)$ we may use Lemma \ref{invo}. If $v(r)=s.v(r-1)$, then $\dim(v(r-1))=\dim v(r)$. If we proceed downwards along the terms of the subexpression we
might have a sequence of steps in which either the $B$-orbit does not
change or it changes through the $W$-action, but we will eventually
reach a step at which $v(j)=m(s')(v(j-1))$ with $\dim(v(j))>\dim
v(j-1)$, where we can apply Lemma \ref{invo}. Hence there is a
$B$-orbit $v'$ in the sequence with $\dim v'=\dim v$, and $\phi(v')$
is an involution $w$. Therefore we may reduce to the case in which $v'=v(r-1)$
and $v=v(r)=s.v(r-1)$. The analysis of the decompositions into
$B$-orbits of $P_sv$ shows that we are in
case IIIa. Then $v\in Bsv'sB\subset BsBwBsB$.      

If $\ell(sws)=\ell(w)+2$ then $v\in BswsB$ and we have the statement.  If $sw>w$ and
$\ell(sws)=\ell(w)$ then $sw=ws$ and $v\in BswB\cup BswsB$ so
$\phi(v)$ is an involution. If $sw<w$ and $\ell(w)=\ell(sws)$ then
$sw=ws$ and $v\in BwB\cup BswB$ so $\phi(v)$ is an involution. Let us
assume that $\ell(sws)=\ell(w)-2$ with $s=s_\alpha$. It follows
from the proof of Lemma \ref{le} in this case that
$\phi(m(s)v)=\phi(m(s)v')=w=\phi(v')$. By \cite[Lemma 1.6]{hum} we have 
$w\alpha\in-\Phi^+$ so we may choose representatives $x=u_x\w \v_x$ and
$y=u_y\w \v_y$ of the same
$B$-orbit $v'$ with $u_x, u_y\in U^{w}$, $\v_x, \v_y\in U$ for which $\alpha$
does not occur in $u_x$ and $\v_y$. 

Conjugation by $\s_\alpha$ maps $x$ in $BswB\cup BswsB$, hence
$\s_\alpha x\s_\alpha^{-1}\in v$ and $\phi(v)\in\{sw,sws\}$. On the
other hand, conjugation by $\s_\alpha$ maps $y$ in $BswsB\cup BwsB$, hence
$\s_\alpha y\s_\alpha^{-1}\in v$ and
$\phi(v)\in\{sw,sws\}\cap\{ws,sws\}$ is an involution. Thus we have the
statement for spherical semisimple conjugacy classes.
  
Let us consider the spherical conjugacy class of an element $x\in G$ with
Jordan decomposition $su$. The proof will follow by induction once we
show that the image through $\phi$ of a closed
$B$-orbit is an involution. As in the proof of Lemma \ref{uno}, if
$y=gxg^{-1}$ is a representative of a closed $B$-orbit then $(G_x\cap
g^{-1}Bg)^\circ$ is a Borel subgroup
of $G_x^\circ$ and $u\in Z(G_x^\circ)$ by \cite[\S 3.15]{131}. Thus
the unipotent part $gug^{-1}$ of $y$ lies in $B$ and
$\phi(B.y)=\phi(B.{gsg^{-1}})$. The conjugacy class $\O_s$ is
spherical because it can be identified with
$G/G_s$ that is a quotient of $G/G_s\cap G_u=G/G_x$. By the
first part of the proof $\phi(B.y)$ is an involution. \hfill$\Box$

\begin{remark}{\rm The reader is referred to \cite[\S 1.4, Remark 4]{ccc} for a
    different proof, in characteristic zero, that the image through $\phi$ of the
    dense $B$-orbit $v_0$, denoted by $z(\O)$, is an involution. In the same
    paper, the image of $\phi(w_0)$ of all spherical conjugacy classes
of a simple algebraic group over ${\mathbb C}$ is explicitely computed.} 
\end{remark}

\section{Stationary points}

In this Section we shall analyze those elements in
$v\in \V$ for which $\phi(m(s)v)=\phi(v)$ for all simple reflections
$s\in W$. 

\begin{definition}Let $v\in \V$, let $w=\phi(v)$ and let
  $\alpha$ be a simple root. We say that
  $v$ is a {\em stationary point with respect to $\alpha$} if
  $\phi(v)=\phi(m(s_{\alpha})v)$. We say that $v$ is a {\em
  stationary point} if it is a stationary point with respect to all
  simple roots.
\end{definition}
 
Stationary points different from the dense $B$-orbit do not exist in
symmetric conjugacy classes by the results in \cite{RS} but they do exist
in unipotent spherical conjugacy classes:

\begin{example}{\rm Let $\O_{\rm min}$ be the minimal nontrivial unipotent
    conjugacy class in a group $G$ of semisimple rank at least
    $2$. It is well-known that $\O_{\rm min}$ is spherical. If $\beta$ denotes the highest root
    in $\Phi$ then $B.x_{\beta}(1)=X_{\beta}\setminus \{1\}$ is a
    stationary point. Indeed, if $\alpha$ is a simple root $P_{\alpha}=Bs_{\alpha}X_{\alpha}\cup B$ and 
$P_{\alpha}.x_{\beta}(t)\subset B.X_{s_{\alpha}(\beta)}\cup
    B.x_{\beta}(1)\subset B$ so $\phi(B.x_{\beta}(1))=\phi(m(s_{\alpha})(B.x_{\beta}(1))=1$.}
\end{example}

The following lemmas describe stationary points with respect to
a simple root. 

\begin{lemma}Let $v\in \V$ with $w=\phi(v)$. Let $\alpha$ be a simple root
  such that $s_{\alpha}w<w$ in the Bruhat order. Then
 $v$ is  a stationary point with respect to $\alpha$.
\end{lemma}
\proof Let us put $s=s_\alpha$. If $sw<w$ then $\phi(m(s)v)\in\{sws,sw,ws,w\}$ and
$\phi(m(s)v)\ge \phi(v)$. If $sw=ws$ the statement follows because
$ws<w$ and $sws=w$. Otherwise it follows because $sw, ws,sws< w$.\hfill$\Box$
\smallskip

\begin{lemma}\label{cnes-stationary}Let $v\in \V$ with $w=\phi(v)$. Let $\alpha$ be a simple root
  such that $s_{\alpha}w>w$ in the Bruhat order. Let $x=u\w\v\in v$ with
  $u\in U^w$, $t\in T$, and  $\v\in U$. Then $v$ is  a stationary point with respect
  to $\alpha$ if and only if the following conditions hold:
\begin{enumerate}
\item  $s_{\alpha}w=ws_{\alpha}$;
\item $\v\in P^u_\alpha$, the unipotent radical of $P_\alpha$;
\item $X_{\pm \alpha}$ commutes with
  $\w\in N(T)$.
\end{enumerate}
\end{lemma}
\proof Let $v\in \V$ be a stationary point with respect to
$\alpha\in\Delta$, let $x\in v$ with unique decomposition $x=u\w\v$
and let $s_\alpha w>w$. It follows from Theorem \ref{image} that $w$ is an involution. 
We have either $s_\alpha w
s_{\alpha}>s_\alpha w$ or $s_\alpha w s_\alpha= w$. If the first case
were possible, we would have
$\s_\alpha x \s_{\alpha}^{-1}\in Bs_\alpha ws_\alpha B$
and $v$ would not be a stationary point. Hence 1 holds,
$w\alpha=\alpha$ and $\alpha$ does not occur in $u\in U^w$. 

Let us consider $y=\s_{\alpha} x \s_{\alpha}^{-1}$. The element 
$$y=(\s_{\alpha}u\s_{\alpha}^{-1})
(\s_{\alpha}\w\s_{\alpha}^{-1})(\s_{\alpha}\v\s_{\alpha}^{-1})\in B
w (\s_{\alpha}\v\s_{\alpha}^{-1}).$$
If $\alpha$ would occur in $\v$ then by \cite[Lemmas 8.1.4, 8.3.7]{springer}
we would have $y\in Bws_\alpha B\cap P_\alpha v$ with $s_\alpha
w>\phi(m(s_\alpha)v)$, a contradiction. Hence 2 holds for
any representative $x\in v$.    

Let then $x\in v$, let $l\in k\setminus\{0\}$ and let $x_1=x_\alpha(l)x
x_{\alpha}(-l)=u_1\w_1\v_1$. Since $\alpha\in \Delta$ and $u,\v\in
P_\alpha^u$ we have $\w_1=\w$ and $\v_1=(\w^{-1}x_\alpha(l)\w)\v x_{\alpha}(-l)$. The root $\alpha$ may
not occur in $\v_1$. By Chevalley commutator formula this is possible only if
$\w^{-1}x_\alpha(l)\w=x_{\alpha}(l)$, that is, only if 3 holds for $X_\alpha$. Let $l\in k\setminus\{0\}$ and let $x_2=x_{-\alpha}(l)x
x_{-\alpha}(-l)$. Since $\alpha\in \Delta$ and $u,\v\in
P_\alpha^u$ we have $x_{-\alpha}(l)u
x_{-\alpha}(-l)\in U$ so
$$x_2\in U\w (\w^{-1}x_{-\alpha}(l)\w)\v x_{-\alpha}(-l).$$
If $(\w^{-1}x_{-\alpha}(l)\w)\v x_{-\alpha}(-l)\not\in U$ we would
have $\phi(B.x_2)\geq ws_\alpha\geq w$ hence $(\w^{-1}x_{-\alpha}(l)\w)\v x_{-\alpha}(-l)\in U$.
By Chevalley commutator formula this is possible only if
$\w^{-1}x_{-\alpha}(l)\w=x_{-\alpha}(l)$, that is, only if 3 holds for
$X_{-\alpha}$.

\smallskip

Let $x$ satisfy  1, 2, and 3. Then $P_\alpha v=P_\alpha
x\subset B\s_\alpha X_\alpha.x\cup v$.  Properties 2 and 3 imply that
$B\s_\alpha X_\alpha.x\subset BwB$ and $v$ is stationary.\hfill$\Box$

\smallskip

\begin{lemma}Let $\O$ be a spherical conjugacy class, let $v\in \V$ be a stationary point and let  $w=\phi(v)$. Let
  $\Pi=\{\alpha\in\Delta~|~w(\alpha)=\alpha\}$ and
$w_\Pi$ be the longest element in
  $W_\Pi$. Then $w=w_\Pi w_0$.
\end{lemma}
\proof By Theorem \ref{image} the
element $w\in W$ is an involution and by Lemma \ref{cnes-stationary}
if $w\alpha\in\Phi^+$ then
$w\alpha=\alpha$ for every $\alpha\in\Delta$. The statement follows
from  \cite[Proposition  3.5]{results}.\hfill$\Box$

\begin{example}{\rm Let $G$ be of type $A_n$ and let $\O$ be a
    spherical conjugacy class. Then the image through $\phi$ of the
    dense $B$-orbit is 
    $w=w_0w_\Pi$ for some $\Pi\subset\Delta$. The set $\Pi$ must be stabilized by
    $-w_0$ because for $\alpha\in\Pi$ we have 
$$\alpha=w\alpha=w_0w_\Pi\alpha\in -w_0\Pi.$$ Besides, if $\alpha_{j}$ lies in
    $\Pi$ then $\alpha_j=w(\alpha_j)=-w_\Pi(\alpha_{n-j+1})$ so
    $\alpha_j$ and $\alpha_{n-j}$ must lie in the same connected
    component of $\Pi$. Hence,
    $\Pi=\{\alpha_t,\alpha_{t+1},\cdots,\alpha_{n-t+1}\}$ for some
    $t$ and 
$$w=(s_{\beta_t}\cdots s_{\beta_{[\frac{n}{2}]}})w_0=s_{\beta_1}\cdots s_{\beta_{t-1}}$$
where $\beta_1,\cdots,\beta_{[\frac{n}{2}]}$ is the sequence given by
the highest root, the highest root of the root system orthogonal to
$\beta_1$, and so further.}
\end{example}

The example above shows that, for any stationary point, $\Pi$ has
to be invariant with respect to $-w_0$. Besides, the restriction of
$w_0$ to $\Phi(\Pi)$ always coincides with $w_\Pi$.  

If $w=\phi(v)$ for some stationary point $v\in \V$,  
the involution $w$ may be written as a product of
reflections with respect to mutually orthogonal roots
$\gamma_1,\ldots,\gamma_m$ so $U_w$ (notation as in Section
\ref{preliminaries}) is the subgroup generated by the 
root subgroups $X_\gamma$ with $(\gamma,\gamma_j)=0$ for every $j$.
In other words, $U_w$ is the subgroup generated by $X_\gamma$ for
$\gamma\in\Phi(\Pi)$ and $U^w$ is normalized by $U_w$. 

\begin{lemma}\label{centra}Let $\O=\O_x$ be a spherical conjugacy
  class in $G$. Let $v\in \V$ be a stationary point. If $x=u\w \v\in U^wN(T)U$ then $B_x\subset u T^wU_w u^{-1}$. 
\end{lemma}
\pf Let $g=\v^w\v_wt\in B_x$ with $\v^w\in U^w$, $\v_w\in U_w$ and $t\in T$, respectively.
Then, 
$$\v^w\v_wt u\w \v=\v^w(\v_w t u t^{-1}\v_w^{-1})\w
(\w^{-1}\v_w\w)(\w^{-1}t\w)\v$$
$$=u\w \v \v^w\v_wt\in U^w\w B$$

By normality of $U^w$ and uniqueness of the decomposition we have
\begin{itemize}
\item $\v^w=u(\v_w t u^{-1} t^{-1}\v_w^{-1})$ so  $\v^w$ is completely
determined by $\v_w$, $t$ and $u$; 
\item $(\w^{-1}\v_w\w)(\w^{-1}t\w)\v= \v \v^w\v_wt$.
\end{itemize}

In particular $g=u\v_wtu^{-1}$ hence $B_x\subseteq
uU_wTu^{-1}$. Uniqueness of the  decomposition in $B=TU$ implies that
$\w^{-1}t\w=t$ so that $t\in T^w$ and $B_x\subseteq u T^wU_w u^{-1}$. 
\hfill$\Box$

\section{The dense $B$-orbit}

We shall turn our attention to the special stationary point given by the
dense $B$-orbit $v_0$. 

We will first analyze the possible $\Pi$ for which $\phi(v_0)=w_0w_\Pi$. We already know that they
are subsets of $\Delta$ for which the restriction of $w_0$ to
$\Phi(\Pi)$ coincides with the longest element of the parabolic
subgroup $W_\Pi$ of $W$. Next step
will be to show which connected components of $\Pi$ may not consist of
isolated roots.

\begin{lemma}\label{quali-no}Let $\O$ be a spherical conjugacy class,
  let $v_0$ be its dense $B$-orbit and let $w=w_0w_\Pi=\phi(v_0)$. 
There is no connected component in $\Pi$ consisting of an
isolated root $\alpha\in\Delta$ such that there is $\beta\in\Delta$
with the following properties:
\begin{itemize}
\item $\beta$ has the same length as $\alpha$;
\item $w_0(\beta)=-\beta$; 
\item $\beta\not\perp\alpha$; 
\item $\beta\perp\alpha'$ for every
  $\alpha'\in\Pi\setminus\{\alpha\}$. 
\end{itemize}
\end{lemma}
\pf Let us choose a representative of $v_0$ of the form $x=\w
\v$ and let us assume that there are simple roots $\alpha$ and $\beta$
with the above properties.  
The element $w$ is the product of $w_0s_\alpha$ with the longest
element $w_{\Pi'}$ of the parabolic subgroup $W_{\Pi'}$ of $W$
associated with the complement $\Pi'$ of $\alpha$ in $\Pi$.
We claim that $\v\in P^u_\beta$. Conversely, given
a representative $\s_\beta$ of $s_\beta$ in $N(T)$, we consider
$y=\s_\beta\w \v\s_{\beta}^{-1}$. Then, as $s_\beta$
commutes with $w_0$ and $w_{\Pi'}$ by assumption, we would have, for
some $l\in k\setminus\{0\}$:
$$y=\w_0\s_{\alpha+\beta}\w_{\Pi'}\v_1
x_{-\beta}(l)\v_2\in Bw_0s_{\alpha+\beta} w_{\Pi'}Bs_\beta B;\quad
\v_i\in P_\beta^u.$$
Besides $w_0s_{\alpha+\beta} w_{\Pi'}(\beta)$ 
is positive and different from $\beta$. Thus, $w_0s_{\alpha+\beta}
w_{\Pi'}s_\beta>w_0s_{\alpha+\beta}
w_{\Pi'}$ and in this case $\phi(B.y)=w_0s_{\alpha+\beta}
w_{\Pi'}s_\beta$ would not be an
involution contradicting Lemma \ref{image}. Hence, $\v\in P^u_\beta$ and $\v'=\s_\beta \v\s_{\beta}^{-1}\in U$. 

Let us consider a representative $\s_\alpha$ of $s_\alpha$ in $N(T)$ and
the element:
$$z=\s_\alpha\s_\beta\w \v\s_{\beta}^{-1}\s_\alpha^{-1}=(\w_0\s_\beta\w_{\Pi'})(\s_\alpha \s_\beta
\v\s_\beta^{-1}\s_\alpha^{-1})\in B w_0s_\beta w_{\Pi'}B(\s_\alpha
\v'\s_\alpha^{-1})B.$$ 
Here we have used that $w_0(\alpha)=-\alpha$ because $\alpha$ is an
isolated root in $\Pi$. If $\v'\in P^u_\alpha$ then
$z\in B w_0s_\beta w_{\Pi'}B$. If $\v'\not\in P^u_\alpha$ then 
$$z\in B w_0s_\beta w_{\Pi'}B(\s_\alpha
\v'\s_\alpha^{-1})B\subset B w_0s_\beta w_{\Pi'}B\cup
B w_0s_\beta s_\alpha w_{\Pi'}B.$$
It follows from Lemma \ref{image}
that also in this case $\phi(B.z)= w_0s_\beta w_{\Pi'}$ because $(w_0s_\beta s_\alpha
w_{\Pi'})^2=s_\alpha s_\beta\neq1$. On the other hand,
$\ell(\phi(B.z))=\ell(w_0s_\alpha w_{\Pi'})=\phi(v_0)$ with $\phi(B.z)\not=\phi(v_0)$,
contradicting $\phi(B.z)\leq\phi(v_0)$.\hfill$\Box$ 

\smallskip
 
\begin{corollary}Let $\O$ be a noncentral spherical conjugacy class,
  let $v_0$ be its dense $B$-orbit and let
  $w=w_0w_\Pi=\phi(v_0)$. Then $\Pi$ is either empty or it is one of
  the following subsets of $\Delta$:

\noindent{Type $A_n$}
$$
\circ\cdots\circ--\bullet\cdots\bullet--\circ\cdots\circ\quad\Pi=\{\alpha_l,\cdots,\alpha_{n-l+1}\},\quad
2\leq l\leq\left[\frac{n}{2}\right]$$  

\noindent{Type $B_n$}
$$
\circ\cdots\circ--\bullet\cdots\bullet =>=\bullet\quad\Pi=\{\alpha_l,\cdots,\alpha_{n}\},\quad2\leq l\leq n$$

$$
\bullet--\circ\cdots\bullet--\circ--\bullet\cdots\bullet =>=\bullet
$$
$$\Pi=\{\alpha_1,\alpha_3,\cdots,\alpha_{2l-1},\alpha_{2l+1},\alpha_{2l+2},\cdots,\alpha_{n}\},\quad
1\leq l\leq \frac{n}{2}$$

\noindent{Type $C_n$}
$$
\circ\cdots\circ--\bullet\cdots\bullet
=<=\bullet\quad\Pi_1=\{\alpha_l,\cdots,\alpha_{n}\},\quad 2\leq l\leq n$$

$$
\bullet--\circ\cdots\bullet--\circ--\bullet\cdots\bullet =<=\bullet
$$ 
$$\Pi=\{\alpha_1,\alpha_3,\cdots,\alpha_{2l-1},\alpha_{2l+1},\alpha_{2l+2},\cdots,\alpha_{n}\},\quad
1\leq l\leq \frac{n}{2}$$

\noindent{Type $D_n$}

$$
\begin{array}{ccc}
\circ\cdots\circ--\bullet\cdots\!&\!\bullet &\!--\bullet\\
&|&\\
&\bullet&\\
\end{array}\quad \Pi_1=\{\alpha_{2l+1},\cdots,\alpha_{n}\},\quad 2\leq l\leq\frac{n}{2}-1$$

$$
\begin{array}{ccc}
\bullet--\circ--\bullet--\circ\cdots\circ--\bullet\cdots\!&\!\bullet &\!--\bullet\\
&|&\\
&\bullet&\\
\end{array}
$$
$$\Pi=\{\alpha_1,\alpha_3,\cdots,\alpha_{2l-1},\alpha_{2l+1},\alpha_{2l+2},\cdots,\alpha_{n}\},\quad
1\leq l<\frac{n}{2}$$

\noindent{$D_{2m}$}
$$
\begin{array}{ccc}
\bullet--\circ--\cdots--\bullet--\!&\!\circ &\!--\bullet\quad \Pi=\{\alpha_1,\alpha_3,\cdots,\alpha_{2m-3},\alpha_{2m-1}\}\\
&|&\\
&\circ&\\
\end{array}
$$

$$
\begin{array}{ccc}
\bullet--\circ--\cdots--\bullet--\!&\!\circ &\!--\circ\quad\Pi=\{\alpha_1,\alpha_3,\cdots,\alpha_{2m-3},\alpha_{2m}\}\\
&|&\\
&\bullet&\\
\end{array}
$$

\noindent{$D_{2m+1}$}
$$
\begin{array}{ccc}
\bullet--\circ--\cdots\bullet--\circ--\!&\!\bullet &\!--\circ\quad\Pi=\{\alpha_1,\alpha_3,\cdots,\alpha_{2m-1}\}\\
&|&\\
&\circ&\\
\end{array}
$$

\noindent{Type $E_6$}

$$
\begin{array}{ccc}
\bullet--\bullet--\!&\!\bullet\!&\!--\bullet--\bullet\quad\Pi=\{\alpha_1,\alpha_3,\alpha_4,\alpha_5,\alpha_6\}\\
&|&\\
&\circ&\\
\end{array}
$$

$$
\begin{array}{ccc}
\circ--\bullet--\!&\!\bullet\!&\!--\bullet--\circ\quad\Pi=\{\alpha_3,\alpha_4,\alpha_5\}\\
&|&\\
&\circ&\\
\end{array}
$$

\noindent{Type $E_7$}

$$
\begin{array}{ccc}
\circ--\bullet--\!&\!\bullet\!&\!--\bullet--\bullet--\bullet\quad\Pi=\{\alpha_2,\alpha_3,\alpha_4,\alpha_5,\alpha_6,\alpha_7\}\\
&|&\\
&\bullet&\\
\end{array}
$$

$$
\begin{array}{ccc}
\circ--\bullet--\!&\!\bullet\!&\!--\bullet--\circ--\bullet\quad\Pi=\{\alpha_2,\alpha_3,\alpha_4,\alpha_5,\alpha_7\}\\
&|&\\
&\bullet&\\
\end{array}
$$

$$
\begin{array}{ccc}
\circ--\bullet--\!&\!\bullet\!&\!--\bullet--\circ--\circ\quad\Pi=\{\alpha_2,\alpha_3,\alpha_4,\alpha_5\}\\
&|&\\
&\bullet&\\
\end{array}
$$

$$
\begin{array}{ccc}
\circ--\circ--\!&\!\circ\!&\!--\bullet--\circ--\bullet\quad\Pi=\{\alpha_2,\alpha_5,\alpha_7\}\\
&|&\\
&\bullet&\\
\end{array}
$$

\noindent{Type $E_8$}

$$
\begin{array}{ccc}
\bullet--\bullet--\!&\!\bullet\!&\!--\bullet--\bullet--\bullet--\circ\\
&|&\\
&\bullet&\\
\end{array}
$$
$$\Pi=\{\alpha_1,\alpha_2,\alpha_3,\alpha_4,\alpha_5,\alpha_6,\alpha_7\}$$

$$
\begin{array}{ccc}
\circ--\bullet--\!&\!\bullet\!&\!--\bullet--\bullet--\bullet--\circ\quad\Pi=\{\alpha_2,\alpha_3,\alpha_4,\alpha_5,\alpha_6,\alpha_7\}\\
&|&\\
&\bullet&\\
\end{array}
$$

$$
\begin{array}{ccc}
\circ--\bullet--\!&\!\bullet\!&\!--\bullet--\circ--\circ--\circ\quad\Pi=\{\alpha_2,\alpha_3,\alpha_4,\alpha_5\}\\
&|&\\
&\bullet&\\
\end{array}
$$

\noindent{Type $F_4$}
$$
\bullet--\bullet
=>=\bullet--\circ,\quad\quad\Pi=\{\alpha_1,\alpha_2,\alpha_3\}
$$

$$
\circ--\bullet
=>=\bullet--\bullet,\quad\quad\Pi=\{\alpha_2,\alpha_3,\alpha_4\}
$$

$$
\circ--\bullet =>=\bullet--\circ,\quad\quad\Pi=\{\alpha_2,\alpha_3\}$$

\noindent{Type $G_2$}
$$\bullet\equiv<\equiv\circ,\quad\quad\Pi=\{\alpha_1\}$$

$$\circ\equiv<\equiv\bullet,\quad\quad\Pi=\{\alpha_2\}.$$
\end{corollary}
\pf Most of the restrictions are due to the fact that
$w_\Pi=w_0|_{\Phi(\Pi)}$. 
The isolated roots occurring as connected components of $\Pi$ are
necessarily alternating, or differ by a node, or their length is
different from the length of all adjacent
roots, as $\{\alpha_n\}$ in type $C_n$. \hfill$\Box$

\begin{remark}{\rm All the above diagrams do actually occur when
    $k={\mathbb C}$ (cfr.
    \cite{ccc}) and they are are strictly more than the
    Araki-Satake diagrams for symmetric conjugacy classes (see
    \cite{araki}, \cite{helminck}, \cite[Table 1]{satake}).}
\end{remark}

\begin{theorem}\label{teo}Let $\O$ be a spherical conjugacy class, let
  $v_0$ be its dense $B$-orbit and let
  $w=\phi(v_0)=w_0w_\Pi$. Then 
$$\dim(\O)=\ell(w)+{\rm rk}(1-w).$$
\end{theorem} 
\Pf We have divided the proof into a sequence of lemmas. 

\smallskip


Let us recall that a standard parabolic subgroup
can be naturally attached to $v\in \V$ (\cite[\S 2]{knop-inv}):
$$P(v)=\left\{g\in G~|~g.v=v\right\}.$$
Let $L(v)$ denote its Levi component containing $T$ and let
$\Delta(v)$ be the corresponding subset of $\Delta$: this is the
so-called set of simple roots of $v$. The derived subgroup
$[L(v),L(v)]$ of $L(v)$ fixes a point in $v$ (this is shown in
\cite[Lemma 1(ii)]{closures} and the argument works also in
positive characteristic). We shall relate $\Pi$, $\Delta(v_0)$
and the centralizer of an element $x\in v_0$ in $B$.

\begin{lemma}\label{deltav}Let $\O$ be a spherical conjugacy class, let
  $v_0$ be its dense $B$-orbit and let
  $w=w_0w_\Pi=\phi(v_0)$. Then $\Delta(v_0)\subset \Pi$. 
\end{lemma}
\Proof If a simple root $\alpha$ occurs in $\Delta(v_0)$ then 
there is $x=u\w\v\in v$ for which $X_\alpha\in U\cap G_x$ and by Lemma \ref{centra}
we have $1\neq u^{-1}X_\alpha u\subset X_\alpha U^w\cap U_w$. This is possible
only if $\alpha\in\Pi$ and $u$ normalizes (hence centralizes) $X_\alpha$.
\hfill$\Box$

\smallskip

\begin{remark}{\rm In characteristic zero it follows from \cite[Page
      289]{closures} and \cite[Corollary 3]{pany3} that
      $\Pi=\Delta(v_0)$. We shall show that this is always the
      case.}
\end{remark}

\begin{lemma}\label{inclusion}Let $\O$ be a spherical
  conjugacy class, let $v_0$ be its dense $B$-orbit, let $x\in v_0$ and let
  $w=w_0w_\Pi=\phi(v_0)$. If there exists $\alpha\in\Delta$ such that
 $X_\alpha\subset P^u_\alpha (P_\alpha)_x\cap U$ then
  $\alpha\in\Delta(v_0)$. 
\end{lemma}
\Pf If for every $r\in k$ there exists $u_r\in P^u_\alpha$ such that 
$g_r=x_\alpha(r)u_r\in (P_\alpha)_x$ this would hold for every choice
of $x\in v_0$ so by Lemma
\ref{centra} we have $\alpha\in\Pi$. Besides, the analysis of the $B$-orbits
in $P_\alpha v_0$ shows that we are either in case I or in case IIa
because the subgroup $H=\psi((P_{\alpha})_x)$ would contain
$\psi(X_\alpha)$. Since $v_0$ is the dense
$B$-orbit, we cannot be in case IIa so 
$P_\alpha v_0=v_0$ and $\alpha\in\Delta(v_0)$.\hfill$\Box$

\smallskip

Let us recall that the depth ${\rm dp}(\beta)$ of $\beta\in\Phi^+$ is the minimal length
of a $\sigma\in W$ for which $\sigma(\beta)\in -\Phi^+$
(\cite{bb}). Then ${\rm dp}(\beta)-1$ is the minimal length
of a $\sigma'\in W$ for which $\sigma'(\beta)\in\Delta$. 

\begin{lemma}\label{shape-of-x}Let $\O$ be a spherical conjugacy class, let $v_0\in \V$
  be the dense $B$-orbit, let $\phi(v_0)=w_\Pi w_0$ and let $L(v_0)$ be
  the Levi factor of $P(v_0)$ containing $T$. There exists $x\in v_0$ such that 
\begin{enumerate}
\item $x$  is centralized by $[L(v_0),L(v_0)]$;
\item $x$ has decomposition $x=\w\v\in U^wwB$;
\item no root in $\Phi^+(\Pi)$ occurs in the expression of $\v$.  
\end{enumerate}
\end{lemma}
\Pf Let $x_1\in v_0$ be an element that is centralized by $[L(v_0),
  L(v_0)]$. If $x_1=u_1\w\v_1$ commutation with $X_{\pm\alpha}$ for
  $\alpha\in\Delta(v_0)\subset \Pi$ and Lemma \ref{cnes-stationary}
  imply that $u_1,\,\w$ and $\v_1$ are centralized by $X_{\pm
  \alpha}$, hence replacing $x_1$ by its $B$-conjugate
  $u_1^{-1}xu_1=x=\w \v$ yields the first and the second statement. Let us assume that for a fixed ordering
  of the positive roots the root $\gamma\in\Phi(\Pi)$ occurs in the
  expression of $\v$ and let us
  assume that $\gamma$ is of minimal depth in $\Phi(\Pi)$ with this
  property. By Lemma \ref{cnes-stationary} the root $\gamma$ is not simple. 
Then, there exists $\sigma=s_{i_1}\cdots s_{i_r}\in W_\Pi$ such that
  $\sigma(\gamma)=\alpha\in\Pi$. Minimality of depth implies that for every other root
  $\gamma'\in\Phi(\Pi)$ occurring in $\v$ and every $j\ge 2$, we have
  $s_{i_j}s_{i_{j+1}}\cdots s_{i_r}(\gamma')\in\Phi^+$. Then for every
  representative $\sigmad\in N(T)$ we would have, for some $t\in T$,  
$$x'=\sigmad x\sigmad^{-1}=\w t\sigmad\v\sigmad^{-1}\in B w_\Pi w_0B$$ so
  the $B$-orbit represented by $x'$ would be stationary with
  $\sigmad\v\sigmad^{-1}\not\in P^u_\alpha$ contradicting Lemma
  \ref{cnes-stationary}.\hfill$\Box$

\smallskip

\begin{lemma}\label{pi=delta}Let $\O$ be a spherical conjugacy class, let
  $v_0$ be its dense $B$-orbit and let $w=w_0w_\Pi=\phi(v_0)$.  
Then $\Delta(v_0)$ contains all long simple roots in $\Pi$. 
In particular, if $\Phi$ is simply laced $\Delta(v_0)= \Pi$ and there
is $x\in v_0$ for which $U\cap G_x=U_w$.
\end{lemma} 
\Proof If $\Pi=\emptyset$ there is nothing to prove; if $\Pi=\Delta$ then
  $w=w_0w_\Pi=1=\phi(v_0)$ so $\O$ is central, $P(v_0)=G$ and
  $\Pi=\Delta=\Delta(v_0)$. We shall assume 
  throughout the proof that $\Pi\neq\emptyset,\Delta$.

Let $x=\w\v$ be as in Lemma \ref{shape-of-x}. If $\Pi$ is strictly
larger than $\Delta(v_0)$ by Lemma \ref{inclusion} then there exists
$\alpha\in\Pi\setminus\Delta(v_0)$ such that $X_\alpha$
is not contained in the centralizer of $x$. Since $X_\alpha\subset G_{\w}$ by Lemma
\ref{cnes-stationary} we conclude that $X_\alpha\not\subset G_{\v}$. By Chevalley's commutator
formula this may happen only if there occurs in $\v$ a positive root
$\gamma$ with $\gamma+\alpha\in \Phi$. Indeed, in this case
the structure constants $c_{\alpha,\gamma,1,1}$ in Chevalley's
commutator formula (\cite[Proposition 9.5.3]{springer}) would be nonzero.
By Lemma \ref{shape-of-x} such a
root $\gamma$ does not lie in $\Phi^+(\Pi)$.  

The basic idea of the proof is to show that if the above happens 
then there exists $v\in \V$ such that $\phi(v)$ is not an involution,
contradicting Theorem \ref{image}. The proof consists in the
construction of an element $\sigma=s_{i_t}\cdots s_{i_1}\in W$ such that:

\begin{enumerate}
\item $\sigma(\gamma)=-\alpha_{i_t}$;
\item the roots $\gamma'_j=s_{i_1}\cdots s_{i_j}\alpha_{i_{j+1}}$ for 
$j<t-1$ do not occur in  $\v$;
\item $\sigma w_0w_\Pi\sigma^{-1}(\alpha_{i_t})\in\Phi^+$ and is different from $\alpha_{i_t}$.
\end{enumerate}

The existence of a $\sigma$ satisfying conditions 1 and 2 guarantees
that for any representative $\sigmad$ of
$\sigma$ in $N(T)$ we have:
$$\sigmad x \sigmad^{-1}\in B\sigma w_0w_\Pi\sigma^{-1}BX_{-\alpha_{i_t}}B=B\sigma w_0w_\Pi\sigma^{-1} Bs_{\alpha_{i_t}}B.$$

Condition 3 guarantees that
$\sigma w_0 w_\Pi\sigma^{-1}<\sigma w w_\Pi\sigma^{-1}s_{i_t}$ and $s_{i_t}$ does
not commute with  $\sigma w_0
w_\Pi\sigma^{-1}$. Then $\phi(B.\sigmad x \sigmad^{-1})=\sigma w_0
w_\Pi\sigma^{-1}s_{i_t}$ that is not an involution.

We will deal with the different possibilities for $\Pi$ separately.

\smallskip

If $\Phi$ of type $A_n$ then
$\Pi=\{\alpha_l,\alpha_{l+1},\ldots,\,\alpha_{n-l+1}\}$ for $1\leq l\leq
\left[\frac{n+1}{2}\right]$. If $\alpha_i\not\in\Delta(v_0)$ for some
$i$ with $l\leq i\leq n-l+1$ then there occurs in $\v$ 
either $\mu_t=\sum_{j=t}^{i-1}\alpha_j$ with $t\leq l-1$ or
$\nu_t=\sum_{j=i+1}^t\alpha_j$ with $t\ge n-l+2$. 
Let $\gamma$ be a root of this form occurring in $\v$ and of minimal height. If
$\gamma=\mu_t$ we consider
$\sigma=\prod_{j=t}^{i-1}s_j$. Then $\sigma(\gamma)=-\alpha_t$ so
condition 1 is satisfied. If for some $r$ with $t<r\leq i-1$ the root 
$\gamma_r=(\prod_{j=i-1}^{r-1}s_j)\alpha_{r}$ would occur in $\v$ we
would  contradict minimality of the height of $\gamma$ so condition 2
is satisfied. Besides, if we put $\vartheta=-w_0$ we have:
$$
\sigma  w_0w_\Pi\sigma^{-1}(\alpha_t)=\sigma
w_\Pi\vartheta(\mu_t)=\alpha_{n-t+1}+{\mbox{ terms independent of
}}\alpha_{n-t+1}
$$
so $\sigma$ satisfies  condition 3. The case $\gamma=\nu_t$ is treated
similarly.  

\smallskip

Let $\Phi$ be of type $B_n$. Then $\Pi$ is either $\Pi_1=\{\alpha_l,\ldots,\alpha_n\}$  with $l\leq n$ of type
$B_{n-l+1}$ or $A_1$; or it is the union of $\Pi_1$ with $|\Pi_1|$ of the
same parity as $n$ and alternating isolated simple
roots. The element $w_\Pi$ is either $w_{\Pi_1}$ or the product of the
reflections corresponding to those isolated simple roots with $w_{\Pi_1}$.

Let us assume that there is $\alpha_i\in\Pi_1\setminus\Delta(v_0)$
with $i\leq n-1$. Let $\gamma\in\Phi^+$ occur in 
$\v$  with $\alpha_i+\gamma\in \Phi$. It follows from Lemma
\ref{shape-of-x} that $\gamma$ is either
$\mu_j=\sum_{p=j}^{i-1}\alpha_p$ for some $j<l$ or
$\nu_j=\sum_{p=j}^i\alpha_p+2\sum_{q=i+1}^n\alpha_q$ for some $j<l$. 
Let us consider a $\mu_j$ with $j$ maximal (that is, $\mu_j$ is of
minimal height of this form). Then we can rule it out by using
$\sigma_j=\prod_{p=j}^is_p$. Condition 2 follows from minimality of
the heigth of $\mu_j$ and
Lemma \ref{shape-of-x}, and condition 3 is checked observing that the
coefficient of
$\alpha_n$ in the expression of $\sigma_jw_\Pi(\gamma)$ is always
positive. Let us then consider $\nu_j$ occurring in $\v$ with $j$
maximal. Then we may rule it out by using
$\tau_j=(\prod_{p=j}^ns_p)(\prod_{q=n-1}^{i+1}s_q)$. Condition 2 and 3
follow from maximality of $j$ and 
Lemma \ref{shape-of-x} and from the fact that the coefficient of
$\alpha_n$ in the expression of $\sigma_jw_\Pi(\gamma)$ is always
positive. Thus $\Pi_1\setminus\Delta(v_0)\subset\{\alpha_n\}$.

Let us now assume that $\Pi$ is the union of $\Pi_1$ and alternating
isolated simple roots. Suppose that $\alpha_i$, for $1\leq i<n$ is an
isolated root in $\Pi\setminus\Delta(v_0)$. Then there occurs in $\v$
a $\gamma$ that is either 
$\mu_j$ as above, with $j<i-1$; $\mu_j'=\sum_{p=i+1}^j\alpha_p$ for
$i+1<j\leq n$; $\nu_j$ as above, with $j<i$; or
$\nu_j'=\sum_{p=i+1}^{j-1}\alpha_p+2\sum_{q=j}^n\alpha_q$ for $j\ge i+2$.
The first set of roots is ruled out as in the case in which
$\alpha_i\in\Pi_1$; the second set of roots is ruled out by using
$\sigma'_j=\prod_{p=j}^{i+1}s_p$ and minimality of the height; the
third set of roots is ruled out by $\tau_j$ as for $\alpha_i\in\Pi_1$. Here, for
condition 2 the non-occurrence in $\v$ of the previous sets of roots
is needed. Let $\nu'_j$ occur in $\v$ with $j$ maximal. Then we
consider
$\tau_j'=(\prod_{q=j-1}^{n-1}s_q)(\prod_{p=n}^{i+1}s_p)$.
Condition 2 is ensured by the maximality of $j$ and by the
non-occurrence of the roots previously excluded. Condition 3 holds because the
coefficient of $\alpha_i$ in $\tau_j'w_\Pi(\nu_j')$ is positive. Thus,
all long roots in $\Pi$ lie in $\Delta(v_0)$. 

\smallskip

Let $\Phi$ be of type $C_n$. Then $\Pi$ is
$\Pi_1=\{\alpha_l,\ldots,\,\alpha_n\}$ with $l\leq n$ of type
$C_{n-l+1}$ or $A_1$, or the union of $\Pi_1$ with alternating isolated simple
short roots. If $\alpha_n\in\Pi_1\setminus\Delta(v_0)$ then there would occur
in $\v$ a root $\gamma$ of the form $\mu_j=\sum_{p=j}^{n-1}\alpha_p$ with
$j<l$. We rule this out by using $\sigma_j=\prod_{p=j}^{n-1}s_p$,
hence $\alpha_n\in\Delta(v_0)$. 

\smallskip

Let $\Phi$ be of type $D_n$. Then $\Pi$ is either
$\Pi_1=\{\alpha_l,\ldots,\,\alpha_{n-1},\alpha_n\}$ for $1\leq l\leq
n-1$ with $l$ odd; the union of $\Pi_1$ with alternating isolated
simple roots; the set of all $\alpha_j$ with $j$ odd; or the set
of all $\alpha_j$ with $j$ odd and $j\leq n-3$ or $j=n$. 

If $\alpha_i\in\Pi_1\setminus\Delta(v_0)$ for $i<n-1$ there would
occur in $\v$ a root among the following:
$\gamma_j=\sum_{p=j}^{i-1}\alpha_p$ for $j\leq l-1$;  
$\delta_j=\sum_{p=j}^i\alpha_{p}+2\sum_{q=i+1}^{n-2}\alpha_q+\alpha_{n-1}+\alpha_n$
for $j\leq l-1$. Let us consider $\gamma_j$ of minimal height among
the $\gamma_t$'s occurring in $\v$. Then
$\sigma_j=\prod_{p=j}^{i-1}s_p$ satisfies the necessary conditions for
$\gamma_j$ for every choice of $\Pi\supset\Pi_1$ so there occurs no
root of type $\gamma_j$ in $\v$. Let us the consider the root of type
$\delta_j$ of minimal height occurring in $\v$. Then the Weyl group
element $(\prod_{p=j}^ns_p)(\prod_{q=n-2}^{i+1}s_q)$ satisfies the
necessary conditions in order to rule out $\delta_j$. Therefore
$\alpha_i\in\Delta(v_0)$ for every $l\leq i\leq n-2$. 
If $\alpha_{n-1}\in\Pi\setminus\Delta(v_0)$ there would
occur in $\v$ a root among the following: 
$\gamma_j=\sum_{p=j}^{n-2}\alpha_p$ for $j\leq l-1$;  
or
$\delta'_j=\sum_{p=j}^{n-2}\alpha_{p}+\alpha_n$ for $j<l$. The roots
of type $\gamma_j$ are ruled out as for $\alpha_i$ for $i<n-1$. In
order to rule out the roots of type $\delta'_j$ we consider the
minimal occurring in $\v$ and we apply
$\tau_j=(\prod_{p=j}^{n-2}s_p)s_n$. Condition 3 is ensured when $\Pi_1\subset\Pi$ because $j<l$ and
when $\Pi$ is a set of alternating roots because the coefficient of 
$\alpha_{n-2}$ in $(\prod_{p=j}^{n-2}s_p)s_nw_\Pi(\delta'_j)$ is
always positive. Thus, $\alpha_{n-1}\in\Delta(v_0)$ and
$\alpha_n\in\Delta(v_0)$ by symmetry.   

Let $\alpha_i$ be an isolated root in $\Pi\setminus\Delta(v_0)$ for
 $i\leq n-2$. Then there would occur in $\v$ a root of type $\gamma_j$
 for $j<i$, $\mu_j=\sum_{p=i+1}^j\alpha_p$ for $i+1\leq j\leq n-1$,
 $\mu_n=\sum_{p=i+1}^{n-2}\alpha_p+\alpha_n$, $\mu'_n=s_n\mu_{n-1}$, 
 $\nu_j=\sum_{p=i+1}^{j-1}\alpha_{p}+2\sum_{q=j}^{n-2}\alpha_q+\alpha_{n-1}+\alpha_n$
 for $j\ge i+2$, or $\delta_j$ for $j<i$.
These roots are handled making use, in this order, of: $\sigma_j$ for
 $\gamma_j$; $\omega_j=\prod_{p=j}^{i+1}s_p$ for $\mu_j$ for $j\leq
 n-1$; $\omega_n=s_n(\prod_{p=n-2}^{i+1}s_p)$ for $\mu_n$;
 $\omega'_n=s_n\omega_{n-1}$ for $\mu'_n$; $(\prod_{p=j}^n
 s_p)(\prod_{q=n-2}^{i+1}s_q)$ for $\nu_j$; $\tau_j$ for $\delta_j$.

In all cases the required conditions are easily verified. In
particular, condition 3 holds, for every choice of $\Pi$ for which
$\alpha_i$ is an isolated root, because
either $\alpha_i$, $\alpha_{i-1}$ or $\alpha_{i+1}$ will occur in
$\sigma w_\Pi\gamma$ with positive coefficient, for $\sigma$ any of
the above Weyl group elements. This concludes the proof for $\Phi$ of
type $D_n$.

\smallskip

Let $\Phi$ be of type $E_6$, so that $\Pi$ is either
$\Pi_1=\{\alpha_1,\alpha_3,\alpha_4,\alpha_5,\alpha_6\}$ or $\Pi_2=
\{\alpha_3,\alpha_4,\alpha_5\}$. For both cases 
the analysis in type $A_n$ shows that
 if $\Pi\neq\Delta(v_0)$ then the roots $\gamma$ occurring in $\v$ for
 which $\gamma+\alpha\in\Phi$ for some $\alpha\in\Pi$ may not lie in the
 root subsystem of type $A_5$ generated by $\Pi_1$. 

Let $\Pi=\Pi_1$ and let us assume that $\alpha_6\not\in\Pi$. Then there
occurs $\gamma$ in $\v$ where $\gamma$ is one of the following roots:
$\mu_1=\alpha_2+\alpha_4+\alpha_5$, $\mu_2=s_3\mu_1$, $\mu_3=s_4\mu_2$,
$\mu_4=s_1\mu_2$, $\mu_5=s_4\mu_4$,
$\mu_6=s_3\mu_5$. All these roots can be ruled out by using,
respectively:
$\sigma_1=s_2s_4s_5$, $\sigma_2=\sigma_1s_3$, $\sigma_3=\sigma_2s_4$,
$\sigma_4=\sigma_2s_1$, $\sigma_5=\sigma_4s_4$,
$\sigma_6=\sigma_5s_3$. Condition 3 is always satisfied because
$\alpha_6$ occurs in the expression of $\sigma_i w_\Pi(\mu_i)$ with
positive coefficient for every $i$. Therefore, $\alpha_6$ (and, symmetrically,
$\alpha_1$) lie in $\Delta(v_0)$ and all roots occurring in $\v$ are
orthogonal to them because $[L(v_0),L(v_0)]$ centralizes $x$.  

Let us assume that $\alpha_4\not\in\Delta(v_0)$. The roots $\alpha_2$
and $\alpha_1+\alpha_2+2\alpha_3+2\alpha_4+2\alpha_5+\alpha_6$ are the
only positive roots that do not lie in $\Phi(\Pi)$, are orthogonal to $\alpha_1$ and
$\alpha_6$ and are such that $\gamma+\alpha_4\in\Phi$. The
first can be ruled out with $\sigma=s_2$ while for the second we may
use $\tau=s_2s_4s_3s_5s_1s_3s_6s_5$. Since there are no positive
roots $\gamma$ orthogonal to $\alpha_1,\,\alpha_4$ and $\alpha_6$, that
do not lie in $\Phi(\Pi)$ and for which
$\gamma+\alpha_3$ or $\gamma+\alpha_5\in\Phi$ we conclude that
$\Pi=\Delta(v_0)$ and we have the statement for $\Pi=\Pi_1$.  

Let $\Pi=\Pi_2$ and let us assume that
$\alpha_3\not\in\Delta(v_0)$. Then there would occur in $\v$ a root
$\gamma$ among the following roots: $\nu_1=\alpha_2+\alpha_4$,
$\nu_2=s_5\nu_1$, $\nu_3=s_6\nu_2$,
$\nu_4=\sum_{j=1}^5\alpha_j+\alpha_4$, $\nu_5=s_6\nu_4$,
$\nu_6=s_5\nu_5$. They can be ruled out by using $\sigma_1=s_2s_4$,
$\sigma_2=\sigma_1s_5$, $\sigma_3=\sigma_2s_6$,
$\sigma_4=\sigma_2s_3s_4s_1$, $\sigma_5=\sigma_4s_6$ and
$\sigma_6=\sigma_5 s_5$, respectively. Here condition 3 follows
because the coefficient of $\alpha_5$ in $\sigma_i w_\Pi\nu_i$ is
positive. Thus $\alpha_3$ and, by symmetry, $\alpha_5$ lie in $\Delta(v_0)$ and all roots occurring in
$\v$ are orthogonal to them. If $\alpha_4$ would not lie in
$\Delta(v_0)$ there would occur in $\v$ either $\alpha_2$ or
$\sum_{j=1}^6\alpha_j$. The first root is ruled out through $s_2$
while the second is ruled out through $s_2s_3s_4s_5s_1s_6$ and we have
the statement in type $E_6$.    
 
\smallskip

Let $\Phi$ be of type $E_7$. Then $\Pi$ corresponds either to a subdiagram $\Pi_1$
of type $D_6$ (containing all simple roots but $\alpha_1$), a
subdiagram $\Pi_2$
of type $D_4$ (containing $\alpha_2,\alpha_3,\alpha_4,\alpha_5$), the union of
$\Pi_2$ and $\alpha_7$, or $\Pi_3=\{\alpha_2,\alpha_5,\alpha_7\}$ .

Let $\Pi=\Pi_1$ and let us assume that $\alpha_3\not\in\Delta(v_0)$.
Then there would occur in $\v$ one of the following roots: 
$\alpha_1$; $\gamma_1=\alpha_4+\sum_{i=1}^5\alpha_i$;
$\gamma_2=s_6\gamma_1$; $\gamma_3=s_5\gamma_2$;
$\gamma_4=s_7\gamma_2$; $\gamma_5=s_5\gamma_4$;
$\gamma_6=s_6\gamma_5$;
$\gamma_7=\alpha_1+2\alpha_2+2\alpha_3+4\alpha_4+3\alpha_5+2\alpha_6+\alpha_7$.
None of these roots may actually occur, as it can be shown through
$s_1$; $\sigma_1=s_1s_3s_4s_2s_5s_4$; $\sigma_2=\sigma_1 s_6$;
$\sigma_3=\sigma_2 s_5$; $\sigma_4=\sigma_2 s_7$; $\sigma_5=\sigma_4
s_5$; $\sigma_6=\sigma_5 s_6$ and
$\sigma_7=\sigma_1s_3(\prod_{j=6}^4s_j)s_2(\prod_{p=7}^4s_p)$,
respectively. Condition 2 holds because all simple reflections except
from the first one in the decomposition of $\sigma_j$ lie in
$W_\Pi$. Condition 3 holds for all $\sigma_j$ because the
coefficient of $\alpha_7$ is positive in
$\sigma_jw_\Pi\gamma_j$. Therefore $\alpha_3\in\Delta(v_0)$ and all
roots in $\v$ are orthogonal to it. If $\alpha_2$ would not lie in
$\Delta(v_0)$ there would occur in $\v$ one of the following
roots: $\delta_1=\alpha_1+\alpha_3+\alpha_4$; $\delta_2=s_5\delta_1$;
$\delta_3=s_6\delta_2$; $\delta_4=s_7\delta_3$;
$\delta_5=\alpha_1+\alpha_2+2\alpha_3+3\alpha_4+2\alpha_5+\alpha_6$;
$\delta_6=s_7\delta_5$; $\delta_7=s_6\delta_6$;
$\delta_8=s_5\delta_7$. None of them may occur, as it is shown through
$\tau_1=s_1s_3s_4$; $\tau_2=\tau_1s_5$; $\tau_3=\tau_2s_6$;
$\tau_4=\tau_3s_7$;
$\tau_5=\tau_1s_2(\prod_{j=5}^3s_j)(\prod_{p=6}^4s_p)$. It follows
that $\alpha_2\in\Delta(v_0)$. Since the only positive root that
is orthogonal to $\alpha_2$ and $\alpha_3$ and lies is
$\Phi\setminus\Phi(\Pi_1)$ is the highest root in $\Phi(E_7)$, we have
the statement for $\Pi_1$. 

Let $\Pi$ be either $\Pi_2$ or the union of $\Pi_2$ and
$\alpha_7$. The discussion concerning type $D_6$ shows that the sum of
a root in $\Phi(\Pi_1)$ occurring in $\v$ with a root in $\Pi_2$ is
never a root. If $\alpha_3\not\in\Delta(v_0)$ 
there would occur in $\v$ some of the $\gamma_i$'s that have been
discussed for $\Pi=\Pi_1$. The Weyl group elements $s_1$; $\sigma_i$
for $i=1,2$ still satisfy the necessary conditions. In order to rule
out $\gamma_3$ we may use $s_1s_3\sigma_3$; for $\gamma_4$ we may use
$\sigma_4'=(\prod_{j=4}^3s_j)(\prod_{p=5}^4s_p)(\prod_{j=6}^4s_j)(\prod_{p=7}^5s_p)s_2s_1s_4$;
for $\gamma_5$ we may use
$\sigma'_5=s_3s_1s_4s_2(\prod_{j=5}^4s_j)(\prod_{p=6}^5s_p)(\prod_{j=7}^6s_j)s_4s_5$;
for $\gamma_6$ we may use $\sigma'_6=\sigma'_5s_5s_4$; for $\gamma_7$
we may use $\sigma'_7=\sigma_1s_3s_1(\prod_{j=6}^4s_j)s_2(\prod_{p=7}^4s_p)$. The necessary conditions hold for both
choices of $\Pi$. Thus, $\alpha_3\in \Delta(v_0)$. If
$\alpha_2\not\in\Delta(v_0)$ there would occur in $\v$ one of the
$\delta_j$'s we mentioned when discussing $\Pi=\Pi_1$. We may use
$\tau_i$ for $i=1,2,3$ to exclude the first three roots; for $\tau_4$
condition 2 is no longer satisfied and we use
$\tau'_4=(\prod_{j=7}^3s_j)s_1$; for $\delta_j$ for $j=5,6$ we may
use $\tau'_j=s_2s_4s_3s_1(\prod_{p=j}^3s_p)(\prod_{q=j+1}^4s_q)$; for
$\delta_7$ we use
$\tau'_7=s_2s_4s_3s_1s_5s_4(\prod_{p=6}^3s_p)(\prod_{q=7}^4s_q)$.
Therefore, $\alpha_2\in\Delta(v_0)$. Arguing as we did for $\Pi=\Pi_1$
we conclude that all $X_\alpha$ for $\alpha\in\Pi$ commute with
$\w\v$. 

Let $\Pi=\{\alpha_2,\alpha_5,\alpha_7\}$. It follows from the
discussion concerning $D_6$ that all roots $\gamma$
occurring in $\v$ and such that $\gamma+\alpha\in\Phi$ for some
$\alpha\in\Pi$ may not lie in the root subsystem $\Phi(\Pi_1)$. Let
us assume that $\alpha_7\not\in\Delta(v_0)$. Then there would occur in
$\v$ one of the following roots:
$\mu_1=\alpha_1+\sum_{j=3}^6\alpha_j$,
$\mu_2=s_2\mu_1$, $\mu_3=s_4\mu_2$, $\mu_4=s_3\mu_3$,
$\mu_5=s_5\mu_3$, $\mu_6=s_3\mu_5$, $\mu_7=s_4\mu_6$,
$\mu_8=s_2\mu_7$. None of these roots may occur since we may use,
respectively: $\sigma_1=s_1\prod_{i=3}^6s_i$; $\sigma_2=\sigma_1s_2$;
$\sigma_3=s_1(\prod_{i=3}^5s_i)s_2(\prod_{i=4}^6s_i)$;
$\sigma_4=(\prod_{i=1}^5s_i)s_4(\prod_{i=2}^4s_i)s_6$;
$\sigma_5=(\prod_{i=1}^5s_i)s_4s_6s_5$;
$\sigma_6=(\prod_{i=5}^2s_i)s_4s_5s_1(\prod_{i=3}^6s_i)$;
$\sigma_7=(\prod_{i=1}^6s_i)s_4s_3s_5s_4$; and
$\sigma_8=(\prod_{i=5}^6(s_is_{i-1}))s_3s_4(\prod_{i=1}^6
s_i)s_3s_4s_2$. Then $\alpha_7\in\Delta(v_0)$ so all roots occurring
in $\v$ are orthogonal to it. Let us assume that $\alpha_2\not\in\Delta(v_0)$. The possible positive roots $\nu$
that are orthogonal to $\alpha_7$, do not lie in $\Phi(\Pi_1)$ and
for which $\alpha_2+\nu\in\Phi$ are:
$\nu_1=\alpha_1+\alpha_3+\alpha_4$,
$\nu_2=s_5\nu_1$,
$\nu_3=\alpha_1+\alpha_2+2\alpha_3+3\alpha_4+2\alpha_5+2\alpha_6+\alpha_7$,
$\nu_4=s_5\nu_3$. These roots are excluded by using
$\tau_1=s_1s_3s_4$, $\tau_2=\tau_1s_5$,
$\tau_3=(\prod_{i=4}^1s_i)(\prod_{i=5}^7s_is_{i-1})s_3s_4$ and
$\tau_4=s_3s_1s_4s_2s_5s_4s_3(\prod_{i=6}^7s_is_{i-1})s_4s_5$. Thus
$\alpha_2\in\Delta(v_0)$ and all roots occurring in
$\v$ are orthogonal to $\alpha_2$ and $\alpha_7$. Since there are no
positive roots $\gamma$ in $\Phi\setminus\Phi(\Pi_1)$ such that
$\gamma\perp\alpha_2,\alpha_7$ and $\gamma+\alpha_5\in\Phi$ we
conclude that $\alpha_5\in\Delta(v_0)$ and we have the statement
for $E_7$. 

\smallskip

Let $\Phi$ be of type $E_8$. Then $\Pi$ corresponds either to a
subdiagram $\Pi_0$ of type $E_7$ (all simple roots but $\alpha_8$), to a subdiagram $\Pi_1$ of type 
$D_6$ (all simple roots but $\alpha_1$, $\alpha_8$), a
subdiagram $\Pi_2$
of type $D_4$ (consisting of $\alpha_2,\alpha_3,\alpha_4,\alpha_5$). 

If $\Pi=\Pi_0$ and $\alpha_3\not\in\Pi_0$ there would occur in $\v$ one
of the following roots: $\gamma_1=\sum_{j=4}^8\alpha_j$;
$\gamma_2=s_2\gamma_1$; 
$\gamma_3=2\alpha_1+2\alpha_2+3\alpha_3+5\alpha_4+4\alpha_5+3\alpha_6+2\alpha_7+\alpha_8$;
$\gamma_4=s_2\gamma_3$; 
$\mu_1=\alpha_4+\sum_{j=1}^8\alpha_j$;
$\mu_2=s_5\mu_1$; $\mu_3=s_6\mu_2$; $\mu_4=s_7\mu_3$;
$\mu_5=\alpha_1+2\alpha_2+2\alpha_3+4\alpha_4+3\alpha_5+2\alpha_6+\alpha_7+\alpha_8$,
$\mu_6=s_7\mu_5$; $\mu_7=s_6\mu_6$; $\mu_8=s_5\mu_7$. 

We exclude $\gamma_1$ and $\gamma_2$ through
$\sigma_1=\prod_{j=8}^4s_j$ and
$\sigma_2=\sigma_1s_2$; $\gamma_3$ through
$\sigma_3=(\prod_{p=8}^1s_p)s_4s_5s_3s_4s_2(\prod_{p=6}^3s_p)(\prod_{j=7}^4s_j)s_1$ 
and $\gamma_4$ through $\sigma_4=\sigma_3 s_2$. We rule out  
$\mu_j$ for $1\leq j\leq 4$ through
$\tau_1=\sigma_2s_3s_1s_4$; $\tau_2=\tau_1s_5$; $\tau_3=\tau_2s_6$;
$\tau_4=\tau_3s_7$. We exclude $\mu_j$ for $j=5,6,7,8$ through
$\tau_5=(\prod_{p=7}^1s_p)s_8s_7s_4s_3s_5s_4s_2(\prod_{q=6}^4s_q)$,
$\tau_6=\tau_5s_7$, $\tau_7=\tau_6s_6$ and $\tau_8=\tau_7s_5$. In
order to verify condition 2 we use Lemma \ref{shape-of-x} and the
non-occurrence of the previously excluded roots. Therefore
$\alpha_3\in\Delta(v_0)$ and all roots occurring in $\v$ are
orthogonal to it. If $\alpha_2\not\in\Delta(v_0)$ there would occur in
$\v$ one of the following roots:
$\beta_1=\alpha_1+\sum_{j=3}^8\alpha_j$;
$\beta_2=\alpha_1+\alpha_2+2\alpha_3+3\alpha_4+2\alpha_5+\alpha_6+\alpha_7+\alpha_8$;
$\beta_3=s_6\beta_2$; $\beta_4=s_7\beta_3$; $\beta_5=s_5\beta_4$;
$\beta_6=s_6\beta_5$; $\beta_7=s_5\beta_3$;
$\beta_8=\alpha_1+2\alpha_2+3\alpha_3+5\alpha_4+4\alpha_5+3\alpha_6+2\alpha_7+\alpha_8$.
All these roots may be ruled out by using, respectively,
$\omega_1=(\prod_{j=8}^3s_j)s_1$; $\omega_2=\omega_1s_2s_4s_3s_5s_4$;
$\omega_3=\omega_2s_6$; $\omega_4=\omega_3s_7$;
$\omega_5=\omega_4s_5$; $\omega_6=\omega_5s_6$;
$\omega_7=\omega_3s_5$;
$\omega_8=s_3s_1(\prod_{p=4}^6s_ps_{p-1})s_2s_4s_3s_1s_7s_8(\prod_{q=5}^4s_qs_{q+1}s_{q+2})s_2(\prod_{q=4}^5s_ps_{p-1})$. 

Thus, $\alpha_2$ and $\alpha_3$ both lie in $\Delta(v_0)$ and all
roots occurring in $\v$ are orthogonal to them. If
$\alpha_5\not\in\Delta(v_0)$ there would occur in $\v$ one of the
following roots: $\nu_1=\sum_{j=6}^8\alpha_j$;
$\nu_2=\alpha_2+\alpha_3+2\alpha_4+\alpha_5+\nu_1$;
$\nu_3=2\alpha_1+2\alpha_2+3\alpha_3+4\alpha_4+3\alpha_5+3\alpha_6+2\alpha_7+\alpha_8$;
$\nu_4=\alpha_2+\alpha_3+2\alpha_4+\alpha_5+\nu_3$.
They are ruled out with $r_1=\prod_{j=8}^6s_j$;
$r_2=(\prod_{j=8}^2s_j)s_4$;
$r_3=(\prod_{j=8}^1s_j)s_4s_3s_5s_4s_2(\prod_{p=6}^3s_p)(s_1s_7s_6)$; 
$r_4=r_3(s_5s_4s_3s_2s_4)$. Thus all roots in $\v$ are orthogonal to
$\alpha_5$ as well. The possible roots $\gamma$  occurring in $\v$
such that $\alpha_6+\gamma\in\Phi$ are then $\alpha_8+\alpha_7$, that
is ruled out by $s_8s_7$, and
$2\alpha_1+2\alpha_2+3\alpha_3+4\alpha_4+3\alpha_5+2\alpha_6+2\alpha_7+\alpha_8$,
that is ruled out by
$\sigma=(\prod_{j=8}^1s_j)s_4s_3s_5s_4s_2s_6s_7(\prod_{p=5}^3s_p)s_1$. The
roots $\gamma$ that might occur in $\v$ for which
$\alpha_7+\gamma\in\Phi$ are $\alpha_8$, that is ruled out by $s_8$,
and
$2\alpha_1+2\alpha_2+3\alpha_3+4\alpha_4+3\alpha_5+2\alpha_6+\alpha_7+\alpha_8$
which is ruled out by $\sigma s_7$. Thus $\alpha_j\in\Delta(v_0)$ for
$j=2,3,5,6,7$. The only positive root that does not lie in $\Phi(E_7)$
and that is orthogonal to these simple roots is the highest root in
$\Phi$, whence the statement.

Let us assume that $\Pi=\Pi_1$. As a consequence of the discussion
concerning type $E_7$ if $\gamma+\alpha\in\Phi$ with
$\alpha\in\Pi$ then $\gamma$ may not lie in the root subsystem
corresponding to $\Pi_0$. Let us assume that $\alpha_3\not\in\Pi$ so
that a root among the $\gamma_j$, $1\leq j\leq4$ and $\mu_p$ for $1\leq
p\leq 8$ introduced when dealing with $\Pi=\Pi_0$ would occur in $\v$. The Weyl group elements $\sigma_1$ and
$\sigma_2$ exclude $\gamma_1$ and $\gamma_2$. In order to
exclude $\gamma_3$ and $\gamma_4$ we use, respectively,
$\sigma_3'=s_8s_7\sigma_3$ and $\sigma_4'=\sigma_3's_2$. Here, in order
to verify condition 2 we use the non-occurrence of $\gamma_1,\gamma_2$.  
The roots $\mu_j$ for $1\leq j\leq 8 $ are excluded by the same Weyl
group elements as for $\Pi=\Pi_0$, therefore $\alpha_3\in\Delta(v_0)$
and all roots occurring in $\v$ are orthogonal to it. If
$\alpha_2\not\in\Delta(v_0)$ there would occur in $\v$ one of the
$\beta_j$'s that were introduced in the discussion of $\Pi=\Pi_0$. 
The roots $\beta_j$ for $1\leq j\leq 8 $ are excluded by the same Weyl
group elements as for $\Pi=\Pi_0$. If
$\alpha_5\not\in\Delta(v_0)$ there would occur in $\v$ one of the
$\nu_i$ for $i=1,2,3,4$. The first two roots are excluded by the same
Weyl group elements $r_1$ and $r_2$ as in case $\Pi=\Pi_0$, where we use for
condition 2 that all roots occurring in $\v$ are orthogonal to
$\alpha_2$ and $\alpha_3$. In order to exclude $\nu_3$ and $\nu_4$ we use
$r_3'=s_8s_7r_3$ and $r_4'=s_8s_7r_4$, respectively. If
$\alpha_6\not\in\Delta(v_0)$ then the possible roots occurring in $\v$
that added to $\alpha_6$ provide a root are ruled out by $s_8s_7$, and
$\sigma'=s_8s_7\sigma$, with notation as for $\Pi=\Pi_0$.  Then the roots
$\gamma$ for which $\alpha_7+\gamma\in\Phi$ are those indicated for
$\Pi=\Pi_0$ and are ruled out by $s_8$ and $\sigma s_7$. As there
are no roots $\gamma$ in $\Phi\setminus\Phi(\Pi_0)$ that are
orthogonal to $\alpha_2,\alpha_3,\alpha_5,\alpha_6$ and for
which $\alpha_4+\gamma\in\Phi$ we have the statement for $\Pi=\Pi_1$. 

Let $\Pi=\Pi_2$. It follows from the discussion concerning $\Phi$ of type $E_7$
that if $\gamma$ occurs in $\v$ and
$\gamma+\alpha\in\Phi$ for some
$\alpha\in\Pi$ then $\gamma$ may not lie in the root subsystem
corresponding to $\Pi_0$. If $\alpha_3\not\in\Delta(v_0)$
there would occur in $\v$ a root among the $\gamma_j$'s and the
$\mu_j$'s that we have previously discussed. Then $\gamma_1$ and
$\gamma_2$ are excluded by $\sigma_1$ and $\sigma_2$; $\gamma_3$ is
excluded by $\sigma_3'$. We rule out $\mu_j$ for $j=1,2,3$ with the
same Weyl group elements $\tau_j$ that we used for $\Pi=\Pi_0$.
Then $\mu_4$ is excluded by $\tau'_4=(\prod_{p=7}^3s_ps_{p+1})s_2(\prod_{q=4}^7s_ps_{p-1})s_1$ and we can rule out $\gamma_4$ by 
$s_1(\prod_{p=3}^8s_p)s_2(\prod_{q=4}^7s_q)(\prod_{p=3}^6s_p)s_2s_4s_5s_3s_4s_2s_1$.
The roots $\mu_j$, for $5\leq j\leq8$ are ruled out, respectively, by:\\
$\tau'_5=(\prod_{p=4}^8s_ps_{p-1})s_2s_4s_3s_1s_5s_6
s_3s_4s_5s_2s_4$;\\
$\tau'_6=(\prod_{p=4}^8s_ps_{p-1})s_2s_4s_3s_1s_6s_5s_7s_6s_3s_4s_5s_2s_4$;\\
$\tau'_7=\tau'_6s_6$ and
$\tau'_8=s_1s_3s_4s_2(\prod_{p=5}^3s_p)s_1(\prod_{q=6}^8s_q)\prod_{p=5}^4(s_ps_{p+1})s_2s_4\prod_{q=7}^5s_q$.
Thus $\alpha_3\in\Delta(v_0)$ and all roots occurring in $\v$ are
orthogonal to it. If $\alpha_2\not\in\Delta(v_0)$ there would occur in
$\v$ one of the $\beta_j$'s previously introduced. We exclude the 
$\beta_j$ for $1\leq j\leq 3$ and $j=7$ by the same Weyl
group elements $\omega_j$ as for $\Pi=\Pi_0$ and $\Pi=\Pi_1$. Then we
may exclude $\beta_j$ for $j=4,5,6$ making use of
$\omega'_4=(\prod_{p=3}^8s_p)s_2s_4s_3s_1(\prod_{q=5}^7s_qs_{q-1})s_3s_4$;
$\omega_5'=\omega_4's_5$ and $\omega'_6=\omega'_5s_6$. Finally, we
exclude $\beta_8$ with the same Weyl group element $\omega_8$ that we
used when $\Pi=\Pi_0$. The order in which we exclude roots allows the
control on condition 2. Thus, $\alpha_2$ and
$\alpha_3\in\Delta(v_0)$ and if $\gamma\in\Phi$ occurs in $\v$ and
$\gamma+\alpha_5\in\Phi$ then $\gamma=\nu_i$ for some $i$, with
notation as for $\Pi=\Pi_0$. The first three roots are ruled out by
the Weyl group elements $r_1,\;r_2,\;r'_3$ that we introduced in the
previous cases.The root $\nu_4$ is ruled out by
$r_4''=s_2s_4s_3s_5s_4s_1s_3(\prod_{p=6}^2s_p)s_4s_7s_8
(\prod_{q=6}^3s_q)s_1(\prod_{p=7}^2s_p)s_4$. 
Thus, all roots occurring in $\v$ are
orthogonal to $\alpha_5,\,\alpha_2,\,\alpha_3$. No root $\gamma$ in
$\Phi\setminus\Phi(\Pi_0)$ for which
$\alpha_4+\gamma\in\Phi$ satisfies these conditions, whence the
statement in type $E_8$.

\smallskip

Let $\Phi$ be of type $F_4$. Then $\Pi$ is either of type $B_3$, $C_3$,
or $B_2$. It follows from the discussion in types $B_n$ and $C_n$ that
all roots $\gamma$ occurring in $\v$ and for which $\alpha+\gamma\in\Phi$ for
some long root $\alpha\in\Pi$ may not lie in any root subsystem of
type $C_3$ or $B_3$ containing $\Pi$. 

If $\Pi=\{\alpha_1,\alpha_2,\alpha_3\}$ is of type $B_3$ and $\alpha_1\not\in\Delta(v_0)$ there would occur in $\v$ either $\nu_1=\alpha_2+\alpha_3+\alpha_4$,
$\nu_2=\alpha_3+\nu_1$, or $\nu_3=\alpha_4+\nu_2$. We rule them out
with $\sigma_1=s_4s_3s_2$, $\sigma_2=\sigma_1s_3$,
$\sigma_3=s_2s_3s_2\sigma_1$ so $\alpha_1\in\Delta(v_0)$ and all roots
occurring in $\v$ are orthogonal to it.  
If $\alpha_2\not\in\Delta(v_0)$, there occurs in $\v$ either
$\alpha_3+\alpha_4$ or $\alpha_1+2\alpha_2+4\alpha_3+2\alpha_2$. We
rule them out through $s_4s_3$ and $s_2s_3s_4s_1s_1s_2s_3$,
respectively. Hence, all long roots in $\Pi$ lie in $\Delta(v_0)$. 

If $\Pi=\{\alpha_2,\alpha_3,\alpha_4\}$ and
$\alpha_2\not\in\Delta(v_0)$ there would occur in $\v$ one of the
following roots: $\alpha_1$, $\gamma_1=\alpha_1+\alpha_2+2\alpha_3$,
$\gamma_2=\gamma_1+\alpha_4$, 
$\gamma_3=s_4\gamma_1$,
$\gamma_4=\alpha_1+2\alpha_2+4\alpha_3+2\alpha_4$. These roots are
excluded by $s_1$, $\sigma_1=s_1s_2s_3$, $\sigma_2=s_3s_2s_1s_2s_4s_3$, $\sigma_3=\sigma_1s_4$, $\sigma_4=s_2s_3s_1s_2s_4s_3$
so $\alpha_2\in\Delta(v_0)$. 

If $\Pi=\{\alpha_2,\alpha_3\}$, then it is contained in a subsystem of
type $C_3$ and in a subsystem of type $B_3$ and we may make use of the
discussion of type $B_n$ and $C_n$. If
$\alpha_2\not\in\Delta(v_0)$ there would occur in $\v$ either
$\alpha_1+\alpha_2+2\alpha_3+\alpha_4$ or
$\alpha_1+2\alpha_2+4\alpha_3+2\alpha_4$. The first root is excluded
by $\sigma=s_3s_2s_1s_2s_4s_3$ and the second is excluded by
$s_2\sigma$.  
 
\smallskip

Let $\Phi$ be of type $G_2$. Then $\Pi$ consists of a simple root.  
If $\Pi=\{\alpha_1\}$ there is nothing to prove.
If $\Pi=\{\alpha_2\}$ and $\Delta(v_0)=\emptyset$
there occurs in $\v$
a root $\gamma$ with $\gamma\in\{\alpha_1,
3\alpha_1+\alpha_2\}$. We may then use
$\sigma=s_1;\,s_2s_1$, respectively, to rule out these possibilites.

\smallskip

The last statement follows from the fact that there is $x\in v_0$ that
is centralized by $[L(\Pi),L(\Pi)]$, so that $U\cap G_x$ is the
unipotent subgroup generated by the root subgroups corresponding to roots in $\Pi$. \hfill$\Box$

\medskip

We have shown so far that $\Pi\setminus\Delta(v_0)$ consists at least
of all long roots in $\Pi$. The method used in the proof of Lemma
\ref{pi=delta} might fail for those 
short roots $\alpha\in\Pi$ for which the $\alpha$-string through some
positive root $\gamma$ is of length at least $3$. When its length is
$3$, that is, in the doubly-laced case, the root $\gamma$ might be
orthogonal to $\Pi$ so condition $3$ might fail. In type $G_2$
condition 2 might fail for all $\sigma\in W$. Next Lemma will deal
with these cases replacing the representative $x\in v_0$ if needed.  

\begin{lemma}\label{pi=deltaBCFG}Let $\Phi$ be of type
  $B_n,\,C_n,\,F_4$ or $G_2$, let $\O$ be a spherical conjugacy class, let
  $v_0$ be its dense $B$-orbit and let $w=w_0w_\Pi=\phi(v_0)$.  
Then $\Delta(v_0)= \Pi$ and there is $x\in v_0$ for which $G_x\cap U=U_w$.
\end{lemma}
\pf Let $x=\w\v\in v_0$ and let $\alpha$ be a
short root in $\Pi$. 
We will show that no positive root $\nu$ with $\nu+\alpha\in\Phi$
may occur in $\v$. For those roots $\nu$ with
$\nu-\alpha\not\in\Phi$ we will argue as in Lemma \ref{pi=delta}. For
the other roots we will show that $\nu+\alpha$ may not occur in any
$X_\alpha$-conjugate of $x$. Then Chevalley commutator formula implies
that $\nu$ may not occur in $\v$. 

\smallskip

We shall deal with the different root systems separately and we will
use the terminology and notation introduced in Lemma \ref{pi=delta}.

Let $\Phi$ be of type $B_n$. The only short root in $\Pi$ is
$\alpha=\alpha_n$. Let us assume that
$\alpha_n\not\in\Delta(v_0)$. The root subgroups that do not commute with
$X_{\alpha_n}$ correspond to the roots $\gamma_j=\sum_{p=j}^{n-1}\alpha_p$
 and $\delta_j=\gamma_j+\alpha_n$ for $j<l$. The roots $\gamma_j$ are ruled out by using
 $\sigma_j=\prod_{p=j}^{n-1}s_p$ and the usual minimality argument. 
If $\Pi$ is the union of $\Pi_1$ and alternating isolated long roots 
then for every $j$ there exists $\alpha\in\Delta(v_0)$ that is not
orthogonal to $\delta_j$, so none of them may occur in $\v$. Therefore if
$\alpha_n\not\in\Delta(v_0)$ we necessarily have $\Pi=\Pi_1$ and some root
$\delta_j$ occurs in $\v$. Let $a\in k$ and let us consider
$x_1=x_{\alpha_n}(a)xx_{\alpha_n}(-a)=\w \v_1$. We claim that
$\delta_j+\alpha_n$ may not occur in $\v_1$ for any choice of
$a\in k$ and for any $j$. Indeed, if $j\geq l$ then the argument in the
proof of Lemma \ref{shape-of-x} shows that
$\delta_j+\alpha_n\in\Phi(\Pi)$ may not occur in $\v_1$. Let
$\delta_j+\alpha_n$ be of
minimal height occurring in $\v_1$. Then $\delta_j+\alpha_n$ is ruled
out by using $\sigma_j s_n$. If some $\delta_j$ would occur in $\v$
then there would exist $a\in k$ for which $\delta_j+\alpha_n$ would
occur in $x_{\alpha_n}(a)\v x_{\alpha_n}(-a)=\v_1$, leading to a
contradiction. The statement is proved for $\Phi$ of type
$B_n$.  

Let $\Phi$ be of type $C_n$ and let
$\Pi_1=\{\alpha_l,\ldots,\,\alpha_n\}$ with $l\leq n$. If
$\alpha_i\in\Pi_1\setminus\Delta(v_0)$ then $i<n$ and there occurs
$\gamma$ in $\v$ among the roots: $\mu_{ji}=\sum_{p=j}^{i-1}\alpha_p$
with $j<l$ and
$\nu_{ji}=\sum_{p=j}^i\alpha_p+2\sum_{q=i+1}^{n-1}\alpha_q+\alpha_n$ for
$j<l$. The first set of roots is ruled through
$\sigma_{ji}=\prod_{p=j}^{i-1}s_p$ while the second set of roots is ruled
out through
$\tau_{ji}=(\prod_{p=j}^{n-1}s_p)(\prod_{q=n}^{i+1}s_q)$. Thus,
$\Pi_1\subset\Delta(v_0)$.  

Let $\Pi$ be the union of $\Pi_1$ and alternating isolated simple
roots and let us assume that one of these, say $\alpha_i$, does not
lie in $\Delta(v_0)$. Then there occurs
in $\v$ one of the roots: $\mu_{ji}$ as before, 
for $j<l$; $\mu_{ji}'=\sum_{p=i+1}^{j}\alpha_p$ for $j\ge i+2$; $\delta_{ji}=\sum_{p=i+1}^{j-1}\alpha_p+2\sum_{q=i}^{n-1}\alpha_q+\alpha_n$
for $j\ge i+2$; $\nu_{ji}$
as before for $j\leq i$; $\nu_{ii}-\alpha_i$.

The roots of type $\mu_{ji}$ are ruled out as when $\alpha_i\in\Pi_1$;
the roots of type $\mu_{ji}'$ are ruled out by using
$\sigma_{ji}'=\prod_{p=j}^{i+1}s_p$; the roots of type $\delta_{ji}$ are
ruled out through
$\omega_{ji}=(\prod_{p=j-1}^{n-1}s_p)(\prod_{q=n}^{i+1}s_q)$ where, for
condition 2, we use the non-occurrence of any $\mu_{ji}'$ in $\v$. Then, 
the roots of type $\nu_{ji}$ for $j<i$
are ruled out as when $\alpha_i\in\Pi_1$. Here, in order to verify
condition 2 we use the non-occurrence in $\v$ of the previous sets of
roots. We see that the root $\nu_{ii}-\alpha_i$ may not occur in $\v$ by
using $\prod_{p=n}^{i+1}s_p$. Hence, $\alpha_i\not\in\Delta(v_0)$ if
and only if $\nu_{ii}$ occurs in $\v$. Let $a\in k$ and let
$x_1=x_{\alpha_i}(a)\w x_{\alpha_i}(-a)=\w\v_1$. We claim that
$\nu_{jj}+\alpha_j$ may not occur in $\v_1$ for any $a\in k$ and for
any $j$. Indeed, if $j\geq l$ we may rule it out
as in Lemma \ref{shape-of-x}. If $j<l$ we consider $\nu_{jj}+\alpha_j$
of minimal height and we rule it out by using $\prod_{p=n}^js_p$. This
implies that $\nu_{ii}$ may not occur in $\v$ so the statement holds
for $\Phi$ of type $C_n$. 

Let $\Phi$ be of type $F_4$ and let us assume that
$\Pi=\{\alpha_1,\alpha_2,\alpha_3\}$. By Lemma \ref{pi=delta} the roots occurring in $\v$ are
orthogonal to $\alpha_1$ and $\alpha_2$. If
$\alpha_3\not\in\Delta(v_0)$ there would occur in $\v$ one of
the following roots: $\alpha_4$ and
$\nu=\alpha_1+2\alpha_2+3\alpha_3+2\alpha_4$. The first root is ruled out
by $s_4$. If $\nu$ occurred in $\v$, there would be an
$X_{\alpha_3}$-conjugate $x_1=\w\v_1$ of $x$ for which $\nu+\alpha_3$
occurs in $\v_1$. This is ruled out by using
$\sigma=s_2s_3s_4s_1s_2s_3$. Condition 2 is ensured by Lemma
\ref{shape-of-x} and by Lemma \ref{pi=delta}, where we showed that
$\alpha_2+\alpha_3+\alpha_4$ and $\alpha_2+2\alpha_3+\alpha_4$ do not
occur in $\v$. Condition 1 and 3 are straightforward. Hence 
$\Delta(v_0)=\{\alpha_1,\alpha_2,\alpha_3\}$. 

Let us assume that $\Pi=\{\alpha_2,\alpha_3\}$. The roots occurring in $\v$ have to be
orthogonal to $\alpha_2$. If
$\alpha_3\not\in\Delta(v_0)$ there would occur in $\v$ one of
the following roots: $\alpha_4$,
$\alpha_1+\alpha_2+\alpha_3+\alpha_4$;
$\mu=\alpha_1+\alpha_2+\alpha_3$ or
$\nu$ as in the previous case. The first and the second root are ruled out,
respectively, by $s_4$ and $s_3s_4s_2s_1$. Then $\mu+\alpha_3$ and $\nu+\alpha_3$ may not
occur in any $X_{\alpha_3}$-conjugate of $\v$ because they are ruled
out, respectively,
by $s_1s_2s_3$ and $\sigma$, whence $\Delta(v_0)=\Pi$. 

Let us assume that $\Pi=\{\alpha_2,\alpha_3,\alpha_4\}$ and that
$\alpha_3\not\in\Delta(v_0)$. The roots occurring in $\v$ have to be
orthogonal to $\alpha_2$ and may not lie in $\Phi(\Pi)$. Therefore 
there would occur in $\v$ one of the following roots:
$\alpha_1+\alpha_2+\alpha_3+\alpha_4$;  $\mu$ or
$\nu$ as in the previous case. The first root is ruled out by
$s_3s_2s_1s_4$. The remaining two roots are handled as in the previous
case: $\mu+\alpha_3$ and $\nu+\alpha_3$ may not occur in any
$X_{\alpha_3}$-conjugate of $\v$ because they are 
ruled out by the same Weyl group element as for
$\Pi=\{\alpha_2,\alpha_3\}$. 
Thus, $\alpha_3\in\Delta(v_0)$ and all roots occurring in $\v$ are
orthogonal to it. There are no roots $\gamma$ that are orthogonal to
$\alpha_3,\alpha_2$, do not lie in $\Pi$ and are such that
$\gamma+\alpha_4\in\Phi$, whence the statement for $\Phi$ of type $F_4$.

Let $\Phi$ be of type $G_2$ let $\Pi=\{\alpha_1\}$. 
If $\Delta(v_0)=\emptyset$ there would occur in $\v$
a root $\gamma\in\{\alpha_2, 
\alpha_1+\alpha_2, 2\alpha_1+\alpha_2\}$. The first two roots are
ruled out through $\sigma=s_2;\,s_1s_2$, respectively. Hence
$\alpha_1\not\in\Pi$ if and only if $2\alpha_1+\alpha_2$ occurs in
$\v$. On the other hand, $3\alpha_1+\alpha_2$ may not occur in any
$X_{\alpha_1}$-conjugate of $\v$ because it is ruled out by
$s_2s_1$. This concludes the proof for $\Phi$ of type $G_2$.\hfill$\Box$ 

\smallskip

\begin{lemma}\label{toro-1} Let $G$ be a simple algebraic group and
    suppose that the longest element $w_0\in W$ acts as $-1$ in its
  standard representation. Let $\O$ be a spherical conjugacy class, let
  $v_0$ be its dense $B$-orbit and let
  $w=w_0w_\Pi=\phi(v_0)$. Then there is $x\in v_0$ for which $U\cap G_x=U_w$ and $T_x^\circ=(T^w)^\circ$ so that
$\dim T.x={\rm rk}(1-w)$.
\end{lemma}
\proof Let $\Gamma(\O)$ be the free abelian subgroup of
all characters of $B$ that arise as weights of eigenvectors of $B$ in
the function field $k(\O)$. Let ${\rm r}(\O)$ denote its
rank as a free abelian group. By \cite[Corollary 1, Corollary
  2]{pany3} in characteristic $0$ or \cite[Lemma 2.1]{knop} for
arbitrary characteristic and by Lemmas \ref{pi=delta},
\ref{pi=deltaBCFG} we have
$\dim(\O)=\ell(w)+{\rm r}(\O).$ Since, for a representative $x\in v_0$
we have $\dim(B_x)\leq \dim T^w+\dim U_w$ 
we have 
\begin{equation}\label{prima}{\rm r}(\O)\geq{\rm rk}(1-w)\end{equation}
where the rank of $1-w$ is the minimal number of reflections
needed in the decomposition of $w$ (cfr. \cite[Page 910]{kostant}), that is, the dimension of
the eigenspace corresponding to $-1$.

Besides, by \cite[Lemma 1(iii)]{closures} the set
$\Pi$ is orthogonal to $\Gamma(\O)$, whence 
\begin{equation}\label{seconda}{\rm r}(\O)\leq n-|\Pi|.\end{equation} 

If $w_0=-1$ the restriction of $w_{\Pi}$ to $\Phi(\Pi)$ is also
$-1$. Then the $-1$ eigenspace of $w_\Pi$ has dimension at least
$|\Pi|$ so the dimension of the $-1$ eigenspace of $w=w_0w_{\Pi}$ is
at most $n-|\Pi|$. This is indeed the exact dimension, for if the
dimension of the $1$-eigenspace of $w$ would be greater than $|\Pi|$, there
would be $\chi$ in ${\mathbb Q}\Delta\setminus{\mathbb Q}\Pi$ such that 
$$\chi=\sum_{\alpha\in
  \Pi}c_\alpha\alpha+\sum_{\beta\in\Delta\setminus\Pi}d_\beta\beta=  w_0w_{\Pi}(\chi)=\sum_{\alpha\in
  \Pi}c_\alpha\alpha-\sum_{\beta\in\Delta\setminus\Pi}d_\beta w_{\Pi}(\beta).$$
We may write $w_\Pi$ as a product of reflections with respect to
  orthogonal roots belonging to $\Phi(\Pi)$ so that $w_{\Pi}(\beta)\in
  \beta+{\mathbb Q}(\Pi)$. If $w_0w_{\Pi}(\chi)=\chi$ we would have
$2\sum_{\beta\in\Delta\setminus \Pi}d_\beta\beta\in{\mathbb Q}\Pi$, 
a contradiction. Thus, ${\rm rk}(1-w)=n-|\Pi|$. Combining this with $(\ref{prima})$ and $(\ref{seconda})$ we obtain the
statement.\hfill$\Box$

\smallskip

\begin{lemma}\label{toronon-1} Let $G$ be a simple algebraic group
  with Dynkin diagram of type $A_n$, $D_{2m+1}$, or $E_6$. Let $\O$ be a
  spherical conjugacy class, let 
  $v_0$ be its dense $B$-orbit and let $w=w_0w_\Pi=\phi(v_0)$. Then
  there is $x\in v_0$ for which $U\cap G_x=U_w$ 
  and  $T_x^\circ=(T^w)^\circ$ so that $\dim T.x={\rm rk}(1-w)$. 
\end{lemma}
\proof Let $x=\w \v$ be as in Lemma \ref{shape-of-x}. Since it is
 centralized by $[L(v_0),L(v_0)]$ we see that ${\alpha}^{\vee}(h)\in T_x$ for every
 $\alpha\in\Pi$. Let $\vartheta$ be the
non-trivial automorphism of $\Phi$ that is equal to $-w_0$. Then
$(T^w)^\circ$ is generated by ${\alpha}^{\vee}(h)$ for $\alpha\in\Pi$ and
${\alpha}^{\vee}(h)(\vartheta\alpha)^{\vee}(h^{-1})$ for
$\alpha\in\Delta\setminus\Pi$. If $(T^w)^\circ\not\subset T_{x}$
there would be a $\gamma$ occurring in $\v$ 
with $(\gamma,\alpha)\neq(\gamma,\vartheta\alpha)$ for some
 $\alpha\in\Delta\setminus\Pi$. We will use the same argument as in the proof of Lemma
 \ref{pi=delta} to show that this cannot be the case.

Let $\Phi$ be of type $A_n$. Then
$\Pi=\{\alpha_l,\ldots,\alpha_{n-l+1}\}$ and by Lemma \ref{pi=delta}
we have $\gamma\perp\Pi$. It follows that
$\gamma=\alpha_j+\cdots+\alpha_t$ with $t\neq n-j+1$ and either $j\leq
l-1$ and $t\ge n-l+1$; or
$j<t\leq l-2$;  or $n-l+3\leq j<t$.  
Let $\gamma$ be the root of minimal height of the first type. If
$t>n-j+1$ we consider $\sigma=s_js_{j+1}\cdots s_t$. Condition 1 is satisfied
by construction. If $s_p\cdots s_t\gamma'=\alpha_{p-1}$ for $p-1>j$
then $\gamma'=\alpha_{p-1}+\cdots\alpha_t$ with $n-t+1<j<p$ so
minimality of the height of $\gamma$ forces condition 2 to hold.
Condition 3 is satisfied because $\alpha_{n-t+1}$
appears in the expression of $\sigma w_0w_\Pi(\gamma)$ with positive
coefficient. If $t<n-j+1$ we shall use $\sigma^{-1}$.  
If $j<t\leq l-2$ or $n-l+3\leq j<t$ the same argument for $\gamma$
minimal applies with $\sigma=s_js_{j+1}\cdots s_t$ and the conditions
are easier to verify. This argument holds also if $\Pi=\emptyset$.

\smallskip

Let $\Phi$ be of type $D_{2m+1}$. If $\Pi_1\not=\emptyset$ (notation
as in Lemma \ref{pi=delta}) then any positive root with
$(\gamma,\alpha)\neq(\gamma,\vartheta\alpha)$ is of the form
$\gamma=\sum_{p=j}^{n-2}\alpha_p+\alpha_q$ for $q=2m, 2m+1$ and
none of these roots is orthogonal to $\Pi_1\subset\Pi$. 

If $\Pi_1=\emptyset$ we consider $\gamma$ of minimal height among
the roots $\gamma$ described above occurring in $\v$ 
and $\sigma=s_j\cdots s_{n-2}s_q$. Then conditions 1 and 2
are satisfied by construction and minimality of height and condition
3 follows because $\alpha_q$ occurs with positive coefficient in
$\sigma w_0w_\Pi\gamma$. 

\smallskip

Let $\Phi$ be of type $E_6$. If $\gamma$ is a positive root for which 
$(\gamma,\alpha_1)\neq(\gamma,\alpha_6)$ then either $\gamma$ or
$\vartheta(\gamma)$ lies in
$\{\alpha_1, \alpha_3,
\alpha_1+\alpha_3,\alpha_3+\alpha_4,
\alpha_2+\alpha_3+\alpha_4,\alpha_1+\alpha_3+\alpha_4, 
\alpha_1+\alpha_2+\alpha_3+\alpha_4, 
\alpha_2+\alpha_3+2\alpha_4+2\alpha_5+\alpha_6,
\alpha_1+\alpha_2+\alpha_3+2\alpha_4+2\alpha_5+\alpha_6\}$
We shall rule out these roots, their image through $\vartheta$ can be
handled symmetrically. The roots ${\alpha_i}$ with $i=1,3$ 
may not occur either because $\alpha_i\in\Pi$, or because
$\alpha_i\not\perp\Pi$, or because they are ruled out using $s_i$.
Besides, $\gamma\perp\Pi$ hence 
$\vartheta\gamma\perp\Pi$ because $\Pi$ is $\vartheta$-stable. 
We may rule out the remaining roots by using 
$\sigma=s_3s_1$,$s_4s_3,s_2s_4s_3$ , $s_4s_3s_1$, $s_2s_4s_3s_1$,
$s_3s_4s_2s_5s_4s_6s_5$, $s_3s_4s_2s_5s_4s_6s_5s_1$. The conditions
in the proof of Lemma \ref{pi=delta}
are easily verified. Let us observe that some roots need only be
handled when $\Pi=\emptyset$, that is, when $w=w_0$.
Thus if a positive root $\gamma$ is such that
$(\gamma,\alpha)\neq(\gamma,\vartheta\alpha)$ for some $\alpha\in\Delta$ 
we necessarily have $\alpha=\alpha_3$ or $\alpha=\vartheta(\alpha_3)$ and either
$\gamma$ or 
$\vartheta(\gamma)$ lies in
$\{\gamma_1=\alpha_1+\alpha_3+\alpha_4+\alpha_5,
\gamma_2=\alpha_1+\alpha_2+\alpha_3+\alpha_4+\alpha_5,
\gamma_3=\alpha_2+\alpha_3+2\alpha_4+\alpha_5+\alpha_6\}$.
These roots can only occur for $\Pi=\emptyset$
and in this case we may use, respectively, 
$\sigma_1=s_5s_4s_3s_1$; 
$\sigma_2=s_2s_4s_3s_1s_5$ and
$\sigma_3=s_5s_4s_2s_3s_4s_5s_6$. By symmetry we deduce that no root
with $\gamma\neq\vartheta\gamma$ may occur in
$\v$. This concludes the proof of the Lemma.\hfill$\Box$

\smallbreak

It follows from Lemmas \ref{pi=delta}, \ref{pi=deltaBCFG}, \ref{toro-1} and
\ref{toronon-1} that there is $x\in v_0$ such that
$(B_x)^\circ=(T^w)^\circ U_w$ and this concludes the proof of Theorem \ref{teo}. \hfill$\Box$ 

\smallskip   

\begin{corollary}\label{dense}Let $\O$ be a spherical conjugacy class, let
  $v_0\in\V$ be the dense $B$-orbit and let $w=w_0w_\Pi=\phi(v_0)$. Then
  $v_0=\O\cap BwB$.
\end{corollary}
\Pf The arguments in Lemmas \ref{shape-of-x} depend only
    on the image of a $B$-orbit through $\phi$. If for some $v\in\V$ we have $\phi(v)=w$ there is a
    representative $x=\w\v$ of $v$ that is centralized by 
$[L(v),L(v)]$ with no occurrence of roots in $\Phi(\Pi)$ in
$\v$. The centralizer of $x$ in $U$ will be contained in $U_w$ and $\Delta(v)\subset \Pi$.  
Proceeding as in Lemmas \ref{pi=delta} and \ref{pi=deltaBCFG} we
    conclude that $\Pi=\Delta(v)=\Delta(v_0)$. The arguments in the
    proof of Lemmas \ref{toro-1} and \ref{toronon-1} involve only the
    image of $\phi$, hence  $\dim v=\dim v_0$, whence $v=v_0$.\hfill$\Box$

\smallbreak

\begin{remark}{\rm If $\O$ is a symmetric conjugacy class over an
  algebraically closed field of odd or zero characteristic  
Theorem \ref{teo}  follows from \cite[Proposition 3.9, Theorem 4.6, Theorem
  7.1]{RS} and Corollary \ref{dense} follows from \cite[Theorem 7.11,
  Lemma 7.12, Theorem 7.13]{RS}. 
If $\O$ is a spherical conjugacy class over an
algebraically closed field of characteristic zero Theorem \ref{teo} is
\cite[Theorem 1]{ccc} and Corollary \ref{dense} is \cite[Corollary 26]{ccc}.} 
\end{remark}

\begin{corollary}\label{chi}Let $\O$ be a spherical conjugacy class with dense
  $B$-orbit $v_0$ and $\phi(v_0)=w_0w_\Pi$. Then there is $x=\w\v\in v_0$ for which all roots
  occurring in $\v$ are fixed by $\vartheta$ and orthogonal to $\Pi$.
\end{corollary}

\begin{corollary}\label{pesi2}Let $\lambda$ be a  weight of an irreducible
  representation occurring in the $G$-module decomposition of the ring
  of regular functions on $\O$ over ${\mathbb C}$. Then
  $\lambda^*=-w_0\lambda=\lambda$. In particular, $\Gamma(\O)$ is
  symmetric.
\end{corollary}
\Proof  By \cite[\S 6]{pany3} we have $\Gamma(\O)={\rm Ann}(T_x)$. 
It follows from Corollary \ref{chi} and Lemmas \ref{toro-1},
\ref{toronon-1} that $(\lambda,
\alpha^\vee)=(\lambda,\vartheta\alpha^\vee)$ for
$\alpha\in\Delta\setminus\Pi$ and $0=(\lambda,
\alpha^\vee)=(\lambda,-\vartheta\alpha^\vee)=(\lambda,\vartheta\alpha^\vee)$
for $\alpha\in\Pi$. Hence $\vartheta\lambda=-w_0\lambda=\lambda$. \hfill$\Box$

\smallskip
\begin{corollary}\label{pesi1}Let $\lambda$ be a  weight of an irreducible
  representation occurring in the $G$-module decomposition of the ring
  of regular functions on $\O$ over ${\mathbb C}$. Let $v_0$ be the
  dense $B$-orbit and $w=w_0w_\Pi=\phi(v_0)$. Then $\lambda\in{\rm Ker}(1+w)\cap P^+\cap Q$.
\end{corollary}
\proof The inclusion in $P^+$ is obvious, the inclusion in $Q$ follows
from the fact that the ring of regular functions on $\O$ is a $G_{\rm ad}$-module. 
It follows from Corollary \ref{pesi1} that
$-\lambda=w_0\lambda=ww_\Pi\lambda=w\lambda$ since $\lambda\perp\Pi$.\hfill$\Box$

\end{document}